\documentclass[11pt]{article}

 \usepackage{varioref,graphicx, epstopdf, amsmath,subfigure, amsfonts} 
\usepackage{amsmath,amsfonts}
\usepackage{color,hyperref} 

\catcode`\@=11
\@addtoreset{equation}{section}
 \catcode`\@=12

\newtheorem{Theorem}{Theorem}[section]
\newtheorem{Proposition}[Theorem]{Proposition}
\newtheorem{Lemma}[Theorem]{Lemma}
\newtheorem{Corollary}[Theorem]{Corollary}
\newtheorem{Remark}[Theorem]{Remark}

\def\dist{{\mathsf{d}}}

\def\Proof{\noindent{{\bf Proof. }}}
\def\square{\vbox{ \hrule height .4pt
\hbox{\vrule width .4pt height 7pt \kern 7pt \vrule width .4pt}\hrule height .4pt }}
\def\QED{\hfill {$\square$}\bigskip}

\def\fa{\forall}

\def\norm#1{\| #1\|}

\def\tfrac#1#2{{\textstyle\frac{#1}{#2}}}

\def\a{{\bf a}}

\def\d{{\delta}}
\def\s{{\sigma}}

\def\e{\varepsilon}
\def\f{\varphi_{2}}
\def\F{\varphi}

\def\ow{{\overline w}}

\def\uw{{\underline w}}

\def\CC{{\cal{C}}}
\def\MM{{\cal{M}}}

\def\HH{{\cal{H}}}
\def\PP{{\cal{P}}}

\def\AA{{\cal{A}}}

\def\R{{\mathbb{R}}}
\def\N{{\mathbb{N}}}

\title{Multiplicity of layered solutions for Allen-Cahn systems with symmetric double well potential}

\author{\small
Francesca Alessio$^{1}$ and Piero
Montecchiari$^{2}$}

\begin{document}
\small

\maketitle

\maketitle
\centerline{\small{
 Dipartimento di Ingegneria Industriale e Scienze
Matematiche, 
} 
}
\centerline{
\small{Universit\`a Politecnica delle Marche, Via Brecce Bianche -- I 60131 Ancona,
}}
\centerline{
\small{
e-mail:
$^{1}${\tt f.g.alessio@univpm.it}, $^{2}${\tt p.montecchiari@univpm.it}.
}
}
\vskip2truecm

\noindent{
\small{\bf Abstract.}
We study the existence of solutions $u:\R^{3}\to\R^{2}$ for the 
semilinear elliptic systems 
\begin{equation}\label{eq:abs}
	-\Delta u(x,y,z)+\nabla W(u(x,y,z))=0,
\end{equation}
where 
$W:\R^{2}\to\R$ is a double well  symmetric potential.  We use variational methods to show, under generic non degenerate properties of the set of one dimensional heteroclinic  connections
between the two minima $\a_{\pm}$ of $W$,
that
(\ref{eq:abs}) has infinitely many geometrically distinct solutions $u\in C^{2}(\R^{3},\R^{2})$ which satisfy
$u(x,y,z)\to \a_{\pm}$ as ${x\to\pm\infty}$ uniformly with respect to $(y,z)\in\R^{2}$ and which exhibit dihedral symmetries with respect to the variables $y$ and $z$. We also characterize the asymptotic behaviour of these solutions as $|(y,z)|\to +\infty$.}


\bigskip

\noindent{\scriptsize {\bf Mathematics Subject Classification:}
35J60, 35B05, 35B40, 35J20, 34C37.

\vskip3truecm
\noindent{
{\scriptsize $^{1,2}$ Partially supported by the PRIN2009 grant "Critical Point Theory and Perturbative Methods for Nonlinear Differential Equations"}}
}
\vfill\eject

%%%%%%%%%%%%%%%%%%%%%%%%%% INTRO %%%%%%%%%%%%%%%%%%%%%%%%%%%%%%%%%%%%%%%%%%%%
\section{Introduction}
We consider semilinear elliptic system of the form

\begin{equation}\label{eq:eq0}
	-\Delta u+\nabla W(u)=0,\qquad\hbox{on }\R^{3}
\end{equation}
where $u:\R^{3}\to\R^{2}$ and $W\in\CC^{3}(\R^{2},\R)$ satisfies
\begin{description}

\item[$(W_{1})$] there exist $\a_{\pm}\in\R^{2}$ such that $W(\a_{\pm})=0$, $W(\xi)> 0$ for every $\xi\in\R^{2}\setminus\{\a_{\pm}\}$ and the Hessian matrixes $\nabla^{2}W(\a_{\pm})$ are definite positive;

\item[$(W_{2})$] 
there exists $R>1$ such that $\nabla W(\xi)\cdot\xi>0$ for every $|\xi|>R$, from which $\mu_0=\inf_{|\xi|>R}W(\xi)>0$;

\item[$(W_{3})$] $W(-\xi_{1},\xi_{2})=W(\xi_{1},\xi_{2})=W(\xi_{1},-\xi_{2})$ for all $(\xi_{1},\xi_{2})\in\R^{2}$;
\end{description}
 In the sequel, without loss of generality, we will assume that ${\bf a}_{\pm}=(\pm 1,0)$.
\bigskip

Models of this kind are used in various fields of Physics, Chemistry or Biology to describe the behaviour of two phase systems.
The two components of the function $u$  represent different order parameters and the two phases ${\bf a}_{\pm}$  are energetically favorite equilibria. In particular (\ref{eq:eq0}) enters in the study of phase transitions. Considering the reaction-diffusion parabolic system 
$\partial_tu-\varepsilon^2\Delta u+\nabla W(u)=0$,
a formal analysis (see  \cite{[BronsardReitich]} and \cite{[Sternberg2]}) shows that, as  $\varepsilon\to 0^+$, its solutions 
tend pointwise to the global minima of $W$ and sharp phase interfaces are produced. System (\ref{eq:eq0}) appears as first term in the expansion at any point on the interface and the corresponding limit solutions $u$ of (\ref{eq:eq0}), named two layered transition solutions,  satisfy the asymptotic condition
\begin{equation}\label{eq:ete} 
\lim_{x\to\pm\infty}u(x,y,z)=\a_{\pm}\quad\hbox{uniformly w.r.t. }  (y,z)\in\R^{2}.
\end{equation}

We recall that in the scalar situation,  when $u : \R^{n}\to \R$ and $W:\R\to\R$ is a double well potential, problem (\ref{eq:eq0})--(\ref{eq:ete}) has been intensively investigated in the last years. In particular the Gibbons conjecture,  proved in \cite{[BarlowBassGui]}, \cite{[BerestyckiHamelMonneau]} and \cite{[Farina]}, establishes that  in the scalar case
the whole set  of solutions to (\ref{eq:eq0})--(\ref{eq:ete}), can be  obtained by  translations of
the unique odd solution 
of the one dimensional heteroclinic problem: $- \ddot q(x) +W'(q(x))=0$ for $x\in\R$ and $q(\pm\infty)=\pm 1$ (for the more general De Giorgi conjecture we refer to \cite{[GhossoubGui]}, \cite{[AmbrosioCabre]}, \cite{[Savin]}, \cite{[DelPinoKowalczykWeiCR]}, \cite{[DelPinoKowalczykWei]} and the survey \cite{[FarinaValdinoci1]}). These results completely describe the set of solutions to (\ref{eq:eq0})--(\ref{eq:ete}) in the scalar context, showing that the problem is in fact {\it one dimensional}.\smallskip

As firstly shown by S. Alama, L. Bronsard and C. Gui in \cite{[ABG]}, the one dimensional symmetry of the solutions to (\ref{eq:eq0})--(\ref{eq:ete}) is generically lost when one considers the vectorial setting. What is basic in their analysis is the fact that, differently from the scalar situation, in the vectorial case the one dimensional heteroclinic problem can have more
 geometrically distinct solutions. Assuming $(W_1)$, $(W_{2})$, $(W_3)$ (the symmetry of $W$ is in fact required only in the $\xi_{1}$ variable)  in \cite{[ABG]}  the existence of an entire solution to (\ref{eq:eq0})--(\ref{eq:ete}) on $\R^{2}$, asymptotic as $y\to\pm\infty$ to two different one dimensional solutions is obtained. This is done assuming that the set of minimal one dimensional heteroclinic connections is constituted, modulo translation, by $k\geq 2$ distinct elements. 

The symmetry condition on the potential, which assures more compactness in the problem, was dropped by  M. Schatzman in \cite{[Schatzman]}, where the same kind of solutions are obtained,
assuming that the set of minimal one dimensional heteroclinic connections consists of exactly two distinct elements which are supposed to be non degenerate, i.e. the kernels of the corresponding linearized operators are one dimensional. In \cite{[Schatzman]} it is furthermore shown that these kind of finiteness and non degeneracy conditions on the set of minimal connections is generic with respect to the choice of the potential $W$ satisfying ($W_{1}$) and ($W_{2}$).

A refined version of the result in \cite{[ABG]} is given for symmetric potentials in \cite{[AlikFusco0]} by N. D. Alikakos and G. Fusco. In that paper also examples of potentials $W$ satisfying the discreteness assumptions on the minimal connections are given together with numerical simulations (see also in this direction  \cite{[Alikbetchen]} and \cite{[AlikFusco00]}).

 We finally refer to \cite{[AlBrakeSyst]} where, adapting to the vectorial case an Energy constrained variational argument used in \cite{[AlM3bump]}, \cite{[AlMalmost]}, 
 \cite{[AlMbrake]} \cite{[AlMNLS]},  it is shown that (\ref{eq:eq0})--(\ref{eq:ete}) admits infinitely many planar solutions whenever the set of one dimensional minimal heteroclinic solutions is not connected.  These planar solutions exhibit different behaviour with respect to the variable $y$, being periodic in $y$ or asymptotic as $y\to\pm\infty$ to one dimensional heteroclinic (not necessarily minimal) connections. They are classified by different values of an ''energy'' parameter.
 
 These results tell us that if the set of minimal one dimensional solutions has some discreteness properties, then the problem (\ref{eq:eq0})--(\ref{eq:ete}) admits a wide variety of planar solutions. A natural question is whether (\ref{eq:eq0})--(\ref{eq:ete}) admits solutions depending on more then two variables. Following a strategy already used in \cite{[AlMtran]} for non autonomous equations, aim of the present paper is to show that under suitable discreteness and non degeneracy properties of the set of minimal one dimensional solutions (\ref{eq:eq0})--(\ref{eq:ete}) admits in fact a multiplicity of different three dimensional solutions.\bigskip
 
 To precisely describe our results we introduce some notation. Letting $z_{0}\in \CC^{\infty}(\R,\R^{2})$ be any fixed function (odd in the first component and even in the second one)  such that $z_{0}(x)=\a_{+}$ for all $x\geq 1$, we consider on the space
\[
\hat\HH_1=z_{0}+\{q\in H^{1}(\R)^{2}\, |\,  q_{1}(x)x>0\,\,\forall x\not=0 \hbox{ and }q(-x)=(-q_{1}(x),q_{2}(x))\},
\]
the action functional
$$
\varphi_{1}(q)=\int_{\R}\tfrac{1}{2}|\dot q(x)|^{2}+W(q(x))\,dx.
$$ 
Due to the symmetry of $W$, the critical points of $\varphi_{1}$ on $\hat\HH_{1}$ (endowed with the $H^{1}$ topology) are classical one dimensional heteroclinic solution of our problem.
Denoting 
\[m_{1}=\inf_{\hat\HH_1}\varphi_{1}\quad\hbox{and}\quad\MM_{1}=\{ q\in\hat\HH_{1}\,|\, \F_{1}(q)=m_{1}\}\] we ask that $\MM_{1}$ is finite, does not contain
scalar connections and consists of not degenerate critical point of $\varphi_{1}$. More precisely, we require
 \begin{itemize}
\item[$(*)$] $\MM_{1}$ is finite and $q_{2}(0)\not=0$ for all $q\in\MM_{1}$;
\item[$(**)$]  there exist $\omega^{*}>0$ such that for all $q\in\MM_1$ there results
$$
\varphi_1''(q)h\cdot h=\int_{\R}|\dot h(x)|^{2}+\nabla^{2}W(q(x))h(x)\cdot h(x)\, dx\geq \omega^{*}\|h\|^{2}_{L^{2}},
\quad\forall h\in \hat\HH_{1}.
$$
\end{itemize}

Note that, by ($W_{3}$), if $q\in\MM_{1}$ then $q^{*}(x)=(q_{1}(x),|q_{2}(x)|)\in\MM_{1}$ too. By uniqueness
of the solution of the Cauchy problem, this implies that $q_{2}(x)$ has constant sign on $\R$. Then, the condition  $q_{2}(0)\not=0$ in ($*$) is actually equivalent to ask (as in the assumption (H4) in \cite{[AlikFusco0]}) that any scalar connections $Q(x)=(Q_{1}(x),0)\in\hat\HH_{1}$, which always exists by ($W_{3}$), is not a minimum for $\varphi_{1}$. The assumption
$(**)$ is exactly equivalent to the non degeneracy requirement made in \cite{[Schatzman]}. The arguments contained in \cite{[Schatzman]} can be used and adapted to the present context to show that ($*$) and ($**$) hold generically (with respect to the $C^{2}$ topology) for potential $W$ satisfying ($W_{1}$), ($W_{2}$) and ($W_{3}$).\bigskip

Denoting $\bar q(x)=(q_{1}(x),-q_{2}(x))$ and observing that $\bar q\in\MM_{1}$ for any $q\in\MM_{1}$ we can finally state our main result.

\begin{Theorem}\label{T:main}
If $(W_{1})$, $(W_{2})$, $(W_{3})$, $(*)$ and $(**)$ hold true, there exist infinitely many solutions 
of the problem 
\begin{equation}\label{eq:P3}
	\begin{cases}- \Delta v(x,y,z) +\nabla W(v(x,y,z))=0,&
	(x,y,z)\in\R^{3},\\
	{\displaystyle\lim_{x\to\pm\infty}}|v(x,y,z)-\a_{\pm}|=0 
	&\hbox{uniformly w.r.t. }(y,z)\in\R^{2}.\end{cases}
\end{equation}
 More precisely, for every $j\geq 2$ there exists $q\in\MM_{1}$ and  a solution $v_{j}\in \CC^{2}(\R^{3},\R^{2})$ of (\ref{eq:P3}) such that,
denoting 
$\check{ v}_{j}(x,\rho,\theta)=v_{j}(x,\rho\cos\theta,\rho\sin\theta)$, it satisfies
\begin{itemize}
\item[$(i)$] $\check{ v}_{j}$ is periodic in $\theta$ with period $\frac{2\pi}{j}$,
\item[$(ii)$]  $\displaystyle\lim_{\rho\to+\infty}\check{ v}_{j}(x,\rho,\tfrac\pi 2+\tfrac\pi{j}(\tfrac 12+k))=\begin{cases}q(x),&\hbox{if }k\hbox{ is odd},\\
\bar q(x),&\hbox{if }k\hbox{ is even},\end{cases}$ uniformly
w.r.t. $x\in\R$.
\end{itemize}
\end{Theorem}

By $(i)$ the solution $v_{j}$ is invariant under rotation with respect to the $x$ axes of angles which are multiple of $2\pi/j$ and so exhibits dyedral symmetry with respect to the variables $y$ and $z$. Moreover, by $(ii)$, Theorem \ref{T:main} gives information about the asymptotic behaviour of the function $v_{j}$ along direction orthogonal to the $x$ axes. Indeed the function  $\check{ v}_{j}(x,\rho,\theta)$ is asymptotic as $\rho\to+\infty$ to 
 $q$ or $\bar q$ for $\theta$   equal 
 to $\tfrac\pi 2+\tfrac\pi{j}(\tfrac 12+k)$  with $k\in \{0,\ldots, 2j-1\}$ odd or even.

As in
\cite{[AlCMsaddle]}, \cite{[AlMsadd3]} and expecially in \cite{[AlMtran]} (see also \cite{[GuiSchatzman]}),
the proof of Theorem \ref{T:main} uses variational methods  to study an auxiliary problem. Indeed, given $j\in\N$, $j\geq 2$, setting $\bar z=\tan(\tfrac{\pi}{2j})z$ and
$$
\PP_{j}=\{(x,y,z)\in\R^{3}\,|\, (x,y)\in \R\times (-\bar z, \bar z),\, z\geq 0\},
$$
we look for functions $v\in \CC^{2}(\PP_{j})^{2}$ satisfying the symmetry conditions $v(-x,y,z)=(-v_{1}(x,y,z),v_{2}(x,y,z))$ and 
$v(x,-y,z)=(v_{1}(x,y,z),-v_{2}(x,y,z))$ for $(x,y,z)\in \PP_{j}$ and which, for a certain $q\in\MM_{1}$, solve the auxiliary problem
 \begin{equation}\label{eq:tildeP3}
\begin{cases} -\Delta v(x,y,z)+\nabla W(v(x,y,z))=0,&
(x,y,z)\in \PP_{j},\cr 
{\displaystyle\lim_{z\to +\infty}}v(x,\bar z,z)=q(x)&\hbox{uniformly w.r.t }x\in\R\cr
\partial_{\nu}v(x,y,z)=0 & (x,y,z)\in \partial \PP_{j},\cr
{\displaystyle\lim_{x\to\pm\infty}}|v(x,y,z)-\a_{\pm}|=0, 
	&\hbox{uniformly w.r.t. }(y,z).
\end{cases}
 \end{equation}
 If $v$ solves (\ref{eq:tildeP3}) then as $z\to +\infty$ we have
 $
 v(x,-\bar z,z)=(v_{1}(x,\bar z,z),-v_{2}(x,\bar z,z))\to (q_{1}(x),-q_{2}(x))=\bar q(x)
 $ and the entire solution $v_{j}$ on $\R^{3}$ is obtained from $v$
  by recursive reflections of the prism ${\mathcal P}_{j}$ with respect to its faces.\smallskip
  
 To solve (\ref{eq:tildeP3}) we build up a renormalized variational procedure (see \cite{[Rabinowitz3]}, \cite{[Rabinowitz4]} and the monography \cite{[RabinowitzStredulinskylibro]} for the use of renormalized functionals in other contexts)  which takes into account the informations we have on the lower dimensional problems. Even if the general strategy of the proof is similar to the one used in \cite{[AlMtran]}, the maximum principle, which leads in the scalar situation to ordering properties of the solutions, is lost in the present setting and a different deeper analysis is needed.
The proof of Theorem \ref{T:main} is contained in section 5 and refers to a list of properties of one dimensional and two dimensional solutions of (\ref{eq:eq0})--(\ref{eq:ete}) studied in sections 2, 3 and 4.

\begin{Remark}\label{R:costanti}
{\rm We precise some basic 
consequences of the assumptions $(W_{1})-(W_{3})$,
fixing some constants. For all $\xi\in\R^{2}$, we set 
$$
\chi(\xi)=\min\{|\xi-\a_{-}|, |\xi-\a_{+}|\}.
$$
First we note that,  
since $W\in\CC^{2}(\R^{2},\R)$  and $\nabla^{2}W(\a_{\pm})$ are definite positive, then
\begin{equation}\label{BGS}
\hbox{$\forall\, r>0$ $\exists\,\omega_{r}>0$ such that if $\chi(\xi)\le r$ then }W(\xi)\geq \omega_{r}\chi(\xi)^{2}.
\end{equation} 
Then, since $W(\a_{\pm})=0$, $\nabla W(\a_{\pm})=\bf 0$ and $\nabla^{2}W(\a_{\pm})$ are definite positive, we have that there exists ${\overline\d}\in (0,\frac{1}{8})$ two constants $\ow>\uw>0$ such that if $ \chi(\xi)\le 2\overline\d$ then 
\begin{equation}\label{eq:iper}
2\uw|\eta|^{2}\le  \nabla^{2}W(\xi)\eta\cdot\eta\le 2\ow|\xi|^{2} \hbox{ for all }\eta\in\R^{2},
\end{equation}
and
\begin{equation}\label{eq:stimeW}
\uw\chi(\xi)^{2}\leq W(\xi)\leq \ow\chi(\xi)^{2}
 \hbox{ and }|\nabla W(\xi)|\leq 2\ow\chi(\xi).
\end{equation} 
}\end{Remark}

\section{One dimensional solutions}

In this preliminary section  we recall some well known properties of the one dimensional minimal solutions to (\ref{eq:eq0}) verifying (\ref{eq:ete}), i.e., minimal solution to the problem
 \begin{equation*}\label{eq:P1bis}
	%\begin{cases}
	- \ddot q(x) +\nabla W(q(x))=0,%&
	\hbox{ for }x\in\R%,\\
	\hbox{ and }{\displaystyle\lim_{x\to\pm\infty}}|q(x)-\a_{\pm}|=0,
	%&\end{cases}
      \eqno(P_{1})
\end{equation*}  
and we display some consequences of the assumptions $(*)$ and $(**)$.\smallskip

Fixed a function $z_{0}\in C^{\infty}(\R,\R^{2})$ 
such that $\|z_{0}\|_{L^{\infty}}\le R$, $z_{0}$ odd in the first component and even in the second one and such that $z_{0}(t)=\a_{+}$ for all $t\geq 1$, we consider on the space
$$
\HH_1=z_{0}+H^{1}(\R)^{2},
$$
the functional
$$
\varphi_{1}(q)=\int_{\R}\tfrac{1}{2}|\dot q(x)|^{2}+W(q(x))\,dx.
$$
We will study some properties of the minima of $\varphi_{1}$ on $\HH_1$ and we set
$$
m_{1}=\inf_{\HH_1}\varphi_{1}.
$$ 
Endowing $\HH_1$ with the hilbertian structure induced by the map
$Q:H^{1}(\R)^{2}\to \HH_1$, $Q(z)=z_{0}+z$, it is well known  that $\varphi_{1}\in \CC^{2}(\HH_1)$
and that critical points of $\varphi_{1}$ are classical solutions to ($P_{1}$). Moreover, given any interval $I\subset\R$ we set
$$
\varphi_{1,I}(q)=\int_{I}\tfrac{1}{2}|\dot q(x)|^{2}+W(q(x))\,dx.
$$
Note that, for every $I\subseteq\R$, the functional $\F_{1,I}$ is weakly lower semicontinuos on $\HH_1$.

\medskip

\begin{Remark}\label{R:simm}
{\rm
Given $q=(q_{1},q_{2})\in\R^{2}$ we denote
$\hat q=(-q_{1},q_{2})$. Then, setting
$$
\hat\HH_1=\{q\in\HH_1\,|\, q(x)_{1}x>0\,\,\forall x\not=0 \hbox{ and } 
q(-x)= \hat q(x),\, t\in\R\},
$$ 
as in Remark 2.2 in \cite{[AlBrakeSyst]} (see also \cite{[GuiSchatzman]}), we can prove that for all $q\in\HH_1$ there exists $Q\in\hat \HH_1$ such that $\varphi_1(Q)\le \varphi_1(q)$ and so that
$$
m_{1}=\inf_{\HH_1}\varphi_1=\inf_{\hat\HH_1}\varphi_1.
$$
}\end{Remark}
\begin{Remark}\label{R:simm}{\rm   If $q\in\hat\HH_{1}$,  then $q(x)_{1}\geq 0$ for $x\geq 0$ and $q(x)_{1}\leq 0$ for $x<0$. Hence $\chi(q(x))=|q(x)-\bf a_{+}|$ for $x\geq 0$ and $\chi(q(x))=|q(x)-\bf a_{-}|$ for $x< 0$, so that $|q(x)-z_{0}(x)|^{2}\leq 2\chi(q(x))^{2}+2\chi(z_{0}(x))^{2}$ for any $x\in\R$. Then, by (\ref{BGS}) we derive that there exists a constant $C$, depending on $\|q\|_{L^{\infty}}$, such that 
$$
\|q-z_{0}\|_{L^{2}}^{2}\leq C\int_{\R}W(q)\, dx+2\int_{\R}\chi(z_{0})^{2}dx\leq C\F_{1}(q)+2\int_{\R}\chi(z_{0})^{2}dx.
$$
 }\end{Remark}

Let us fix $\delta_{0}\in(0,\overline\d)$ such that 
\begin{equation}\label{eq:lambda0}
\lambda_{0}:=\sqrt{2\uw}\d_{0}(\overline\d-\d_{0})-\frac{\d_{0}^{2}}2(1+2\ow)>0.
\end{equation}
Moreover note that
if  $q\in H^{1}_{loc}(\R)^{2}$
is such that $W(q(t))\geq \mu$ for all $t\in (\sigma,\tau)\subset\R$, $\mu>0$, then
\begin{equation}\label{eq:stime1dim}
    \F_{1,(\sigma,\tau)}(q)	\geq
	{\textstyle{\frac{1}{ 2(\tau-\sigma)}}}{|q(\tau)-q(\sigma)|}^{2}+
	\mu(\tau-\sigma)
	\geq \sqrt{2 \mu} \
	|q(\tau)-q(\sigma)|.
\end{equation}

\begin{Remark}\label{R:boundLinfty+conc}
{\rm
As consequence of Remark \ref{R:simm} and (\ref{eq:stime1dim}), using $(W_{2})$, as in Lemmas 2.1 and 2.2 in \cite{[AlBrakeSyst]}, we obtain that 
there exists $R_{0},\, C_{0},\,T_{0}>0$ such that if $q\in\hat\HH_1$ and $\F_{1}(q)\le m_{1}+\lambda_{0}$  then 
\begin{itemize}
\item[$(i)$] $\|q \|_{L^{\infty}}\leq R_{0}$ and $\|q-z_{0}\|_{H^{1}}\le C_{0}$;
\item[$(ii)$] $|q(x)-\a_{+}|\le\overline\delta$ for all $x>T_{0}$.
\end{itemize}
}\end{Remark}

Using Remark \ref{R:boundLinfty+conc} we plainly obtain 
\begin{Lemma}\label{L:min}
    Let $(q_{n})$ be a sequence in $\hat\HH_1$ such that $\F_{1}(q_{n})\le m_{1}+\lambda_{0}$ for all $n\in\N$.
    Then, there exists $q\in
    \hat\HH_1$ 
    such that, along a subsequence,
    $q_{n}-q\to 0$ weakly in $H^{1}(\R)^{2}$ and $\F(q)\le\displaystyle\liminf_{n\to +\infty} \F(q_{n})$. 
\end{Lemma}
By Lemma \ref{L:min}, we obtain that
$
\MM_{1}=\{q\in\hat\HH_1\,|\, \F_{1}(q)=m_{1}\}
$
 is not empty. 
 Moreover, as we recalled above,  every $q\in\MM_{1}$ is a classical $C^{2}(\R)$ solution to problem ($P_{1}$). 
It is simple to show that the elements of $\MM_{1}$ are uniformly  exponentially asymptotic to the points $\a_{\pm}$.

\begin{Lemma}\label{L:expSOL1dim} 
For every $q\in\MM_{1}$,  $|q(x)-\a_{+}|\le \bar\d e^{-\sqrt{\underbar w/2}(x-T_{0})}$, for all
$x\geq T_{0}$.
\end{Lemma}
\Proof Let $q\in\MM_{1}$. By Remark \ref{R:boundLinfty+conc} we have $|q(x)-\a_{+}|\le\overline\delta$ for $x>T_{0}$. Setting $\phi(x)=|q(x)-\a_{+}|^{2}-\bar\d^{2}e^{-\sqrt{2\underline w}(x-T_{0})}$ and using (\ref{eq:iper}) we recover  $\ddot \phi(x)\geq 2\underline {w}\phi(x)$, $\phi(T_{0})\leq 0$ and $\lim_{x\to+\infty}\phi(x)=0$ which imply $\phi(x)\leq 0$ for $x\geq T_{0}$.\QED

\noindent Lemma \ref{L:min} establishes that every minimizing sequence for $\varphi_1$ over $\hat\HH_1$ is precompact with respect to the weak $H^{1}(\R)^{2}$ topology. As in Lemma 2.4 in \cite{[AlBrakeSyst]} the result can be improved and we have

  \begin{Lemma}\label{L:compK}
  Let $(q_{n})$ be a minimizing sequence for $\varphi_1$ over $\hat\HH_1$.
    Then, there exists $q\in
    \MM_{1}$ such that, along a subsequence, $\|q_{n}-q\|_{H^{1}}\to 0$ as $n\to +\infty$.
    \end{Lemma}

\noindent In particular, by Lemma \ref{L:compK}, for every $r>0$ there exists
$\nu_{r}>0$ such that
\begin{equation}\label{eq:h0}
    \hbox{if
    }q\in\hat\HH_1\hbox{ and }%\F_{1}(q)\leq m_{1}+\p\hbox{ and }
    \dist_{H^{1}}(q,\MM_{1})\geq r\hbox{
    then }\F_{1}(q)\geq m_{1}+\nu_{r}.
\end{equation}     

\noindent We finally discuss some consequences of the assumptions $(*)$ and $(**)$. Recall first the discreteness assumption $(*)$ on $\MM_{1}$: 

\begin{itemize}
\item[$(*)$] $\MM_{1}$ is finite and $q(0)_{2}\not=0$ for all $q\in\MM_{1}$.
\end{itemize}

\noindent Since $(*)$ requires that $\MM_{1}$ is finite,  we have in particular that 
	\begin{equation}\label{eq:d0bis}
	\min_{p,q\in\MM_{1},p\not=q}\|q-p\|_{L^{2}}= 5d_{0}>0.
	\end{equation}

\noindent The assumption $(**)$ asks that 

\begin{itemize}
	\item[$(**)$]  there exist $\omega^{*}>0$ such that for all $q\in\MM_{1}$ there results
$$
\varphi_{1}''(q)h\cdot h=\int_{\R}|\dot h(x)|^{2}+\nabla^{2}W(q(x))h(x)\cdot h(x)\, dx\geq \omega^{*}\|h\|^{2}_{L^{2}}
\quad\forall h\in \hat\HH_{1}.
$$
\end{itemize}

\noindent As consequence, using the Taylor's Formula and  Lemma \ref{L:compK} we obtain the following coerciveness property of $\F_{1}$ 

\begin{Lemma}\label{L:Fmenoc}
There exists $\nu^{*}\in (0,\lambda_{0})$ such that if $q\in\hat\HH_{1}$ is such that $\F_{1}(q)\leq m_{1}+
\nu^{*}$,  then
$$
\F_{1}(q)-m_{1}\geq \frac{\omega^{*}}{4}
\dist_{L^{2}}(q,\MM_{1})^{2}.
$$
\end{Lemma}
\Proof We set $\overline W=\frac16\max_{|\xi|\leq R_{0}}|D^{3}W(\xi)|$ and let $c_{0}$ be the constant of the immersion of $H^{1}(\R)^{2}$ into $L^{\infty}(\R)^{2}$. By (\ref{eq:h0})  there exists $\nu^{*}\in (0,\lambda_{0})$  such that if $q\in\hat\HH_{1}$
and $\F_{1}(q)\leq m_{1}+
\nu^{*}$ then 
\begin{equation*}\label{eq:norma}
\|q-q_{0}\|_{H^{1}(\R)^{2}}\leq\min\{d_{0},\tfrac{\omega^{*}}{8\overline {W} c_{0}}\},
\end{equation*}
for some $q_{0}\in\MM_{1}$. Since, by Remark \ref{R:boundLinfty+conc}
we have $\|q\|_{L^{\infty}}\leq R_{0}$ and $\|q_{0}\|_{L^{\infty}}\leq R_{0}$, 
by the Taylor Formula and $(**)$, we obtain that 
\begin{align*}
&\F_{1}(q)-m_{1}=\F_{1}(q)-\F_{1}( q_{0})= \F_{1}'(q_0)(q-q_0)+\tfrac12 \F_{1}''(q_0)
(q-q_0)(q-q_0)+\\
&\,\,\,+\int_{\R}(W(q)-W(q_0)-\nabla W(q_0)\cdot(q-q_0)-\tfrac12
\nabla^{2}W(q_0)(q-q_0)\cdot(q-q_0)\,dx\\
&\geq \frac{\omega^{*}}2\| q-q_0\|^{2}_{L^{2}}-\overline W\|q-q_0\|^{3}_{L^{3}}\geq \left(\frac{\omega^{*}}2-\overline W c_{0}\|q-q_0\|_{H^{1}}\right) \| q-q_0\|^{2}_{L^{2}}
\end{align*}
and the Lemma follows.
\QED

\section{Two dimensional heteroclinic type solutions}

In this section  we display some results concerning the solutions of  the two dimensional problem
 \begin{equation*}
	\begin{cases}- \Delta v(x,y) +\nabla W(v(x,y))=0,&
	(x,y)\in\R^{2}\\
	{\displaystyle\lim_{x\to\pm\infty}}|v(x,y)-\a_{\pm}|=0, 
	&\hbox{uniformly w.r.t. }y \in \R,\cr
	\end{cases}\eqno{(P_{2})}
\end{equation*}
 which are asymptotic as $y\to\pm\infty$ to $\MM_1$. 

 \medskip
 
\noindent Let us consider the renormalized functional
\begin{align*}
\varphi_{2}(v)&= \int_{\R}(\int_{\R}\tfrac{1}{ 2}|\nabla v(x,y)|^{2}+W(v(x,y))\,dx-m_{1})\, dy\\
	&=\int_{\R}\tfrac12\|\partial_{y}v(\cdot,y)\|_{L^{2}}^{2}+\varphi_1(v(\cdot,y))-m_{1}\, dy
\end{align*}
which is well defined on the space 
$$
\HH_{2}=\{v\in H^{1}_{loc}(\R^{2})^{2}\,|\, v(\cdot,y)\in\hat\HH_1\hbox{ for almost every $y\in\R$}\}.
$$
On $\HH_{2}$, given any interval $I\subset\R$, we will also consider  the functional
\[
\varphi_{2,I}(v)= \int_{I}\tfrac12\|\partial_{y}v(\cdot,y)\|_{L^{2}}^{2}+\varphi_1(v(\cdot,y))-m_{1}\, dy, \quad v\in\HH_{2}.
\]
Note that $\varphi_{2.I}(v)\geq 0$ for all $v\in \HH_{2}$, $I\subseteq\R$ and
if $q\in \MM_1$, then the function $v(x,y)= q(x)$
belongs to $\HH_{2}$ and
$\f(v)=0$, i.e., the  minimal solutions of
($P_1$) are global minima of $\f$ on $\HH_{2}$. \medskip

\noindent We will look for bidimensional solutions of ($P_{2}$) as minima of $\f$ on  suitable subspaces of $\HH_{2}$. 
We recall ({see e.g. \cite{[AlJMprimo]}}) that
$\f$ (and $\varphi_{2,I}$ for $I\subseteq\R$) is weakly lower semicontinuous on $\HH_{2}$ with respect to the $H^{1}_{loc}(\R^{2})^{2}$. Concerning the coerciveness of $\f$, we firstly list some basic estimates.\smallskip

\noindent First we note that if $v\in \HH_{2}$ then
$\varphi_1(v(\cdot,y))\geq m_{1}$ for almost every $y\in\R$ and so
\begin{equation}\label{eq:distorsione}
	\norm{\partial_{y}v}^{2}_{L^{2}}\leq 2\f(v)
	\qquad\fa v\in \HH_{2}.
\end{equation}
Moreover, since $W(\xi)\geq 0$ for any $\xi\in\R^{2}$, 
we have
\begin{align*}
	\varphi_{2,(y_{0},y_{1})}(v)&\geq\int_{y_{0}}^{y_{1}}\F_{1}(v(\cdot,y))-m_{1}\, dy+\tfrac{1}{2}
	\norm{\partial_{y}v}^2_{L^{2}(\R\times(y_{0},y_{1}))}\\
	&
	\geq
	\tfrac{1}{2}\norm{\partial_{x}v}_{L^{2}(\R\times(y_{0},y_{1}))^{2}}^{2} +\tfrac{1}{2}
	\norm{\partial_{y}v}^2_{L^{2}(\R\times(y_{0},y_{1}))} -m_{1}(y_{1}-y_{0})
\end{align*}
from which we derive
\begin{equation}\label{eq:limitatezza}
	\norm{\nabla v}^{2}_{L^{2}(\R\times(y_{0},y_{1}))^{2}}\leq
	2\varphi_{2,(y_{0},y_{1})}(v)+2m_{1}(y_{1}-y_{0})\quad\fa\, v\in \HH_{2},\ (y_{0},y_{1})\subset\R.
\end{equation}
Finally, if $v\in \HH_{2}$ then
$v(x,\cdot)\in H^{1}_{loc}(\R)^{2}$ for almost every $x\in\R$. Therefore,
if $y_{0}<y_{1}\in\R$ then
$v(x,y_{1})-v(x,y_{0})=\int_{y_{0}}^{y_{1}}\partial_{y}v(x,y)\, dy$
holds for  almost every $x\in\R$, whenever $v\in \HH_{2}$. Hence, if $v\in \HH_{2}$ and
 $y_{0}<y_{1}\in\R$, by (\ref{eq:distorsione}) we obtain 
\begin{equation}\label{eq:contd}
	\|v(\cdot,y_{1})-v(\cdot, y_{0})\|_{L^{2}}^{2}	\le
	|y_{1}-y_{0}|\int_{\R}\int_{y_{0}}^{y_{1}}|\partial_{y}v|^{2}\,	dy\,
	dx \leq 2\varphi_{2,(y_{0},y_{1})}(v) |y_{1}-y_{0}|.
\end{equation}

\noindent 
By (\ref{eq:h0}), if $(y_{0},y_{1})\subset\R$ and $v\in
\HH_{2}$ are such that 
$\dist_{H^{1}}(v(\cdot,y),\MM_1)\geq d>0$
for almost every $y\in (y_{0},y_{1})$, then there exists $\nu_{d}>0$ such that $\F_{1}(v(\cdot,y))\geq m_{1}+\nu_{d}$ and hence
\begin{align}
	\nonumber	\varphi_{2,(y_{0},y_{1})}(v)	 &\geq \tfrac{1}{ 2(y_{1}-y_{0})}\int_{\R}(\int_{y_{0}}^{y_{1}}
	|\partial_{y}v(x,y)|\, dy)^{2}\, dx+\nu_{d}(y_{1}-y_{0})\\ \label{eq:stime2dim}	
			&\geq\tfrac{1}{2(y_{1}-y_{0})}
	\|v(\cdot,y_{1})-v(\cdot,y_{0})\|_{L^{2}}^{2}+\nu_{d}(y_{1}-y_{0})\\ \nonumber
	&
	\geq \sqrt{2\nu_{d}}\|v(\cdot,y_{1})-v(\cdot,y_{0})\|_{L^{2}}.
\end{align}

\noindent As consequence of (\ref{eq:contd}) and (\ref{eq:stime2dim}) we get information on the asymptotic behavior  as $y\to\pm\infty$
of the functions in the sublevels of $\f$. More
precisely, as in  Lemma 3.2 in \cite{[AlBrakeSyst]} and  Lemma 3.1 in \cite{[AlMtran]}, we have

 \begin{Lemma}\label{L:asymptot2}
	If $v\in \HH_{2}$ and $\f(v)<+\infty$, then
	there exists $q\in\MM_1$ such that $\|v(\cdot,y)-q\|_{L^{2}}\to 0$ as $y\to +\infty$.
\end{Lemma}

\noindent Considering the symmetry of $W$, we look for solutions $v$ of ($P_{2}$) which satisfy the symmetry condition
$(v(x,-y)_{1},v(x,-y)_{2})=(v(x,y)_{1},-v(x,y)_{2})$ and which connect different elements  of $\MM_1$ as $y\to \pm\infty$.
Hence we define 
$$
\tilde \HH_{2}=\{ v\in\HH_{2}\,|\, v(x,-y)=\tilde v(x,y),\, (x,y)\in\R^{2}\}
$$
where we denote $\tilde v=(v_{1},-v_{2})$ for every $v=(v_{1},v_{2})\in\R^{2}$.\\
We remark that by definition, if $v\in\tilde\HH_{2}$, then, a.e. in $\R^{2}$ there results
\begin{align*}
&v(-x,y)_{1}=-v(x,y)_{1}\hbox{ and }v(-x,y)_{2}=v(x,y)_{2}\hbox{, while }\\
&v(x,-y)_{1}=v(x,y)_{1}\hbox{ and }v(x,-y)_{2}=-v(x,y)_{2},
\end{align*}
from which in particular there results $v(-x,-y)=-v(x,y)$. Hence,  if $v\in\tilde \HH_{2}$ is such that  $\|v(\cdot,y)-q\|_{L^{2}}\to 0$ as $y\to +\infty$ for some $q\in\MM_1$, then $\|v(\cdot,y)-\bar q\|_{L^{2}}\to 0$ as $y\to -\infty$,  where $\bar q(x)=-q(-x)$.\\ 
\begin{Remark}\label{R:l2estimate}
{\rm Using the symmetry properties of the functions in $\tilde\HH_{2}$ we gain  coerciveness property for $\f$. Precisely, we have that there exists a constant $\tilde C>0$ such that
if $(y_{0},y_{1})\subset\R$ and if $v\in\tilde\HH_{2}$ is such that $\|v\|_{L^{\infty}}\leq R$ then
\begin{equation}\label{eq:stimal2}
\|v-z_{0}\|^{2}_{L^{2}(\R\times(y_{0},y_{1}))}\leq 
\tilde C\varphi_{2,(y_{0},y_{1})}(v)+2\tilde C (y_{1}-y_{0})m_{1}+4(y_{1}-y_{0})\int_{\R}\chi(z_{0})^{2}\, dx.
\end{equation}
Indeed, arguing as in Remark \ref{R:simm}, we have that 
 $|v(x,y)-z_{0}(x)|^{2}\leq 2\chi(v(x,y))^{2}+2\chi(z_{0}(x))^{2}$ on $\R^{2}$. By (\ref{BGS}), since
$\|v\|_{L^{\infty}}\leq R$, there exists a constant $\tilde C$ such that $|v(x,y)-z_{0}(x)|^{2}\leq \tilde CW(v(x,y))+2\chi(z_{0}(x))^{2}$ and (\ref{eq:stimal2}) plainly follows. 	}\end{Remark}

Then, we have
\begin{Lemma}\label{L:compattezza}
	Let $(v_{n})$ be a sequence in $\tilde\HH_{2}$ such that
	 $\f(v_{n})\leq C$ and $\|v_{n}\|_{L^{\infty}}\leq R$ for all $n\in\N$ and some $C>0$. Then, 
	 there exists $v\in \tilde\HH_{2}$ such that, up to a subsequence,
	$v_{n}-v\to 0$ weakly 
in $H^{1}(\R\times(-L,L))^{2}$ for every $L>0$. 
\end{Lemma}
\Proof
By (\ref{eq:limitatezza}) and (\ref{eq:stimal2}),  $(v_{n}-z_{0})$ is bounded in
$H^{1}(\R\times(-L,L))^{2}$ for all $L>0$ and so there exists a subsequence
$(v_{n_{k}})$ of $(v_{n})$ 
and a function $v$ such that $v-z_{0}\in\cap_{L>0}H^{1}(\R\times(-L,L))^{2}$ and
$v_{n_{k}}-z_{0}\to v-q_{0}$ weakly $H^{1}(\R\times(-L,L))^{2}$ for every 
$L>0$ and a.e. in $\R^{2}$. 
Since $v-z_{0}\in\cap_{L>0}H^{1}(\R\times(-L,L))^{2}$
we have that $v(\cdot,y)-z_{0}\in H^{1}(\R)^{2}$ for almost every $y\in\R$ and by pointwise convergence, $v\in \tilde\HH_{2}$ follows.\QED

\begin{Remark}\label{R:semicondist}	{\rm Standard semicontinuity arguments (see e.g. the proof of Lemma 3.3 in \cite{[AlMtran]})  show that
if  $(v_{n})_{n\in\N}\subset \tilde\HH_{2}$, $v\in\tilde\HH_{2}$ are such that $v_{n}-v\to 0$ weakly in $H^{1}(\R\times(-L,L))^{2}$ for all $L>0$ and if
$q\in\HH_{1}$ then
$\|v(\cdot,y)-q\|_{L^{2}}\le\displaystyle\liminf_{n\to +\infty}\|v_{n}(\cdot,y)-q\|_{L^{2}}$ for almost every $y\in\R$.}\end{Remark}
\begin{Remark}\label{R:ysimm}
{\rm Note that as consequence of the symmetry properties of the functions in $\tilde\HH_{2}$,we have that if $v\in \tilde\HH_{2}$ then $\dist_{L^{2}}(v(\cdot,{0}),\MM_1)\geq {d_{0}}$.\\
 Indeed, given any $q\in \MM_1$, considering $\bar q(x)=- q(-x)$,  we have $\bar q\in\MM_1$ and then, by $(*)$, $q\not=\bar q$ and  $\|q-\bar q\|_{L^{2}}\geq 2d_{0}$.
 By symmetry, if $v\in\tilde\HH_{2}$, we have $v(-x,0)_{1}=-v(x,0)_{1}$ and $v(x,0)_{2}=0$ and then
 \begin{align*}
 \|
 v(\cdot,0)-\bar q\|_{L^{2}}^{2}
 &=\int_{\R}|v(x,0)_{1}-\bar q(x)_{1}|^{2} +|\bar q(x)_{2}|^{2}dx
 \\&=\int_{\R}|v(x,0)_{1}-q(x)_{1}|^{2} +|q(x)_{2}|^{2}dx= \|v(\cdot,0)-q\|_{L^{2}}^{2},
 \end{align*}
 from which we deduce $\|v(\cdot,0)-q\|_{L^{2}}\geq \frac12\|q-\bar q\|_{L^{2}}\geq {d_{0}}$.
Since $q\in\MM_1$ is arbitrary, we conclude that $\dist_{L^{2}}(v(\cdot,0),\MM_1)\geq {d_{0}}$.
}
\end{Remark}

\noindent Having in mind Lemma  \ref{L:asymptot2}, given $q\in\MM_1$ we look for minima of $\f$ over the class
\begin{align*}
	\tilde\HH_{2,q}&=\{ v\in \tilde\HH_{2}\,|\, \lim_{y\to
	+\infty}\|v(\cdot,y)-q\|_{L^{2}}= 0\}
\end{align*}
and we set
$$
m_{2,q}=\inf_{v\in\tilde\HH_{2,q}}\f(v)\quad\hbox{and}\quad \MM_{2,q}=\{v\in\tilde\HH_{2,q}\,|\,\f(v)=m_{2,q}\}.
$$
\begin{Remark}\label{minimi} {\rm One can plainly see that $m_{2,q}<+\infty$ for all $q\in\MM_1$. Moreover, by  $({*})$ and Lemma \ref{L:asymptot2}, we have $m_{2}=\inf_{\tilde\HH_{2}}\f=\min_{q\in\MM_1}m_{2,q}$. 
We denote
\[
\MM_1^{*}=\{q\in\MM_1\,|\, m_{2,q}=m_{2}\}\hbox{ and } 
\MM_{2}=\{v\in\tilde\HH_{2}\,|\, \f(v)=m_{2}\}
\]
noting that $\MM_{2}=\cup_{q\in\MM_1^{*}}\MM_{2,q}$.
Finally, we set 
\begin{equation}\label{eq:lambdatilde}
\tilde\lambda=\begin{cases}{\displaystyle\min \{m_{2,q}-m_{2},\,q\in\MM_1\setminus\MM_1^{*}\} }&\hbox{if $\MM_1\setminus\MM_1^{*}\not=\emptyset$}\cr
1&\hbox{if $\MM_1\setminus\MM_1^{*}=\emptyset$}\end{cases}
\end{equation}
}
\end{Remark}

\noindent Setting  
\begin{equation}\label{eq:lambdabar}
\bar\lambda=\frac18\min\{{\lambda_{0}};\tilde\lambda;d_{0}\sqrt{\nu_{0}}\}, 
\end{equation}
(where $\lambda_{0}$ is given in (\ref{eq:lambda0}), $\tilde\lambda$  in (\ref{eq:lambdatilde}) and $\nu_{0}$ by (\ref{eq:stime2dim}) corresponding to $d=\frac{d_{0}}2$) and letting
$r_{0}\in (0,\min\{\frac{d_{0}}2,\sqrt{2\bar\lambda}\})$
such that
\begin{equation}\label{eq:defr0}
	\sup\{ \varphi_1(q)\, |\,q\in\hat\HH_1,\,\dist_{H^{1}}(q,\MM_1)\leq 2r_{0}\}\leq m_{1}+\bar\lambda,
\end{equation}
we have  the following concentration result
\begin{Lemma}\label{L:dd2dim} If $L>0$, $v\in\tilde\HH_{2}$, 
	$\varphi_{2,(-L,L)}(v)\leq m_{2}+\bar\lambda$
	 and $\|v(\cdot,y_{0})-q\|_{H^{1}}\leq r_{0}$ for some $y_{0}\in( 0,L)$ and $q\in\MM_1$, then	$$q\in\MM_{1}^{*}\hbox{ and }
	\hbox{$\|v(\cdot,y)-q\|_{L^{2}}\leq  d_{0}$ for all
		$y\in (y_{0},L)$ 
		}
		$$
\end{Lemma}
\Proof 
 Let $v,q,L,y_{0}$ as in the statement. 
We define
$$
		 {w}(x,y)=
		 \begin{cases}		 v(x,y)&\hbox{if }x\in\R\hbox{ and }0\le y\le y_{0}\cr
		 v(x,y_{0})(y_{0}+1-y)+q(x)(y-y_{0})&\hbox{if }x\in\R\hbox{ and }y_{0}\leq 
	 	 y\leq y_{0}+1,\cr
		 q(x)&\hbox{if }x\in\R\hbox{ and } y\geq y_{0}+1.\end{cases}
$$
and $w(x,y)=-w(-x,-y)$ for $x\in\R$ and $y<0$, noting that
 $w\in\tilde\HH_{2,q}$ and so  $\f(w)=2\F_{2,(0,+\infty)}(w)\geq
m_{2,q}\geq m_{2}$. Since $\|w(\cdot,y)-q\|_{H^{1}}\le 2r_{0}$ for  $y\in(y_{0},y_{0}+1)$, by (\ref{eq:defr0})
we obtain 
\begin{equation}\label{eq:ppp}
\F_{2,(y_{0},y_{0}+1)}(w)=\int_{y_{0}}^{y_{0}+1}\tfrac{1}{2}
		\int_{\R}|v(x,y_{0})-q(x)|^{2}dx+
		\varphi_1(w(x,y))-m_{1}\, dy\le  
		2\bar\lambda.
\end{equation} 
In particular, by (\ref{eq:ppp}) and using (\ref{eq:lambdabar}), we derive that 
\[\varphi_{2}(w)=\varphi_{2,(-y_{0},y_{0})}(v)+2\varphi_{2,(y_{0},y_{0}+1)}(w)\leq
m_{2}+5\bar\lambda<m_{2}+\tilde\lambda\] 
and so, by the definition of $\tilde\lambda$ in (\ref{eq:lambdatilde}), that $q\in\MM_{1}^{*}$.\\
Moreover, by (\ref{eq:ppp}) again, we obtain that $\tfrac{m_{2}}2\leq \F_{2,(0,+\infty)}(w)	
	\le \F_{2,(0,L)}(v)-\F_{2,(y_{0},L)}(v)+2\bar\lambda,$
from which we derive $\F_{2,(y_{0},L)}(v)\leq 3\bar\lambda$.\\
If we assume by contradiction that there exists $y_{1}\in (y_{0},L)$ such that
$\|v(\cdot,y_{1})-q\|_{L^{2}}>d_{0}$,
then, by
(\ref{eq:contd}) there exists $(y_{0}',y_{1}')\subset (y_{0},y_{1})$
such that $\|v(\cdot,y_{1}')-v(\cdot,y_{0}')\|_{L^{2}} \geq \frac{d_{0}}2$ and $\frac{d_{0}}2\le \|v(\cdot,y)-q\|_{L^{2}}\le d_{0}$ for
any $y\in (y_{1}',y_{0}')$. In
particular, since $\|p-q\|_{L^{2}}> d_{0}$ for all $p\not=q\in\MM_1$, we obtain
 $
\dist_{H^{1}}(v(\cdot,y),\MM_1)\geq \dist_{L^{2}}(v(\cdot,y),\MM_1)=\|v(\cdot,y)-q\|_{L^{2}}\geq
\tfrac{d_{0}}2
$ 	for almost every $y\in (y_{0}',y_{1}')$. Then, by
(\ref{eq:stime2dim}) and (\ref{eq:lambdabar}), we get the contradiction
$
3\bar\lambda\geq\F_{2,(y_{0},L)}(v)\geq \sqrt{2\nu_{0}}\|v(\cdot,y_{1}')-v(\cdot,y_{0}')\|\geq 
\sqrt{2\nu_{0}}\tfrac{d_{0}}2\geq4\bar\lambda$ which completes the proof.\QED

\noindent Using (\ref{eq:h0}) and (\ref{eq:stime2dim}) we  fix $\bar\ell>0$ such that if $I$ is
any real interval with  length $|I|\geq \bar\ell$ 
then
\begin{equation}\label{eq:lbar}
\hbox{if }v\in \tilde\HH_{2}\hbox{ and }\dist_{H^{1}}(v(\cdot,y),\MM_1)>r_{0}\hbox{ for almost every
}y\in I\hbox{ then }\varphi_{2,I}(v)\geq m_{2}+2\bar\lambda.
\end{equation}
As simple consequence of Lemma \ref{L:dd2dim}, we then obtain

\begin{Lemma}\label{L:concentrazione2dL}
	If $L\geq\bar\ell$, $v\in\tilde\HH_{2}$ and $\varphi_{2,(-L,L)}(v)\leq m_{2}+\bar\lambda$
	then there exists $q\in\MM_1^{*}$ such that $\|v(\cdot,y)-q\|_{L^{2}}\le d_{0}$ for  all $y\in( \bar\ell,L)$.
\end{Lemma}
\Proof   Since, by Remark \ref{R:ysimm},
$\dist_{L^{2}}(v(\cdot,{0}),\MM_1)\geq {d_{0}}$ and since
$\varphi_{2,(-L,L)}(v)\leq m_{2}+\bar\lambda$, by definition of $\bar\ell$ there exists 
$y_{0}\in (0, \bar\ell)$ and $q\in\MM_1$ such that 
$\|v(\cdot,y_{0})-q\|_{H^{1}}\le r_{0}$. Then the Lemma follows by Lemma \ref{L:dd2dim}.
\QED

\noindent In particular
\begin{Lemma}\label{L:concentrazione2d}
	If $v\in\tilde\HH_{2}$ and $\f(v)\leq m_{2}+\bar\lambda$
	then there exists $q\in\MM_1^{*}$ such that $\|v(\cdot,y)-q\|_{L^{2}}\le d_{0}$ for  all $y\geq \bar\ell$.
\end{Lemma}

\noindent The results stated in Lemmas \ref{L:asymptot2} and \ref{L:concentrazione2d} imply 
$$
\{v\in\tilde\HH_{2}\,|\, \f(v)\le m_{2}+\bar\lambda\}\subset\cup_{q\in\MM_1^{*}}\{v\in\tilde\HH_{2,q}\,|\, \f(v)\le m_{2}+\bar\lambda\}
$$
and  Lemma \ref{L:concentrazione2d}, Lemma 
\ref{L:compattezza} allow us to use the direct
method of the Calculus of Variation to show that  $\f$
admits a minimum in every class $\tilde\HH_{2,q}$ with $q\in\MM_1^{*}$.

\begin{Proposition}\label{P:sol2dim}
For every $q\in\MM_1^{*}$ there exists $v\in\MM_{2,q}$.
\end{Proposition}
\Proof Fixed $q\in\MM_1^{*}$, let $(v_{n})_{n\in\N}$ be a minimizing sequence for $\varphi_{2}$ in $\tilde\HH_{2,q}$.
Then, by (\ref{eq:d0bis}), the definition of $\tilde\HH_{2,q}$ and Lemma \ref{L:concentrazione2d}, for all $n\in\N$ there results
		\begin{equation}\label{uno}
		\dist_{L^{2}}(v_{n}(\cdot,y),\MM_1)=\|v_{n}(\cdot,y)-q\|_{L^{2}}\leq  d_{0}\
		\hbox{ for }\, y\geq \bar\ell.
	\end{equation}
Now, by Lemma \ref{L:compattezza}, there exists $ v\in\tilde\HH_{2}$ such that, along a subsequence, $v_{n}\to v$ weakly in $H^{1}(\R\times(-L,L))^{2}$ for every $L>0$ and a.e. in $\R^{2}$. By semicontinuity, $\varphi_{2}(v)\le m_{2}$ and  by (\ref{uno}) and Remark \ref{R:semicondist}, 
$\|v(\cdot,y)-q\|_{L^{2}}\leq  d_{0}$ for $y\geq \bar\ell.$
Then,  by Lemma \ref{L:asymptot2} since $\dist_{L^{2}}(q,\MM_{1}\setminus\{q\})\geq 2d_{0}$, we obtain $\|v(\cdot,y)-q\|_{L^{2}}\to 0$ as $y\to +\infty$. Hence $v\in\tilde\HH_{2,q}$ and $\varphi_{2}(v)=m_{2}$ follows. 
\QED

\begin{Remark}\label{R:trapped}
{\rm  By symmetry of $W$, as in Lemma \ref{L:principiovariazionale} below, it can be proved that every $v\in\MM_{2}$ is a weak solutions of $-\Delta v+\nabla W(v)=0$ on $\R^{2}$.  By (\ref{eq:limitatezza}) and (\ref{eq:stimal2}), using bootstrap arguments, we can conclude that $v$ is a $\CC^{2}$ solutions  which satisfies the symmetry conditions
$$
\hbox{$v(-x,y)=\hat v(x,y)$ and $v(x,-y)=\tilde v(x,y)$ for all $(x,y)\in\R^{2}$.}
$$
Moreover, given $v\in \tilde\HH_{2}$, let $P(v)$ be defined as  
$$
P(v)(x,y)=\begin{cases} \frac{Rv(x,y)}{|v(x,y)|}&\hbox{if $|v(x,y)|>R$,}\cr
 v(x,y)&\hbox{if  $|v(x,y)|\leq R$.}
 \end{cases}
 $$
 The assumption $(W_{2})$ guarantees that $\f(P(v))<
\f(v)$ whenever $v\in \tilde\HH_{2}$ is such that $\|v\|_{L^{\infty}}>R$, hence, since $P(v)\in\tilde\HH_{2}$ for every $v\in \tilde\HH_{2}$, we derive that
$\|v\|_{L^{\infty}}\leq R$ for every $v\in\MM_{2}$.
Hence, by (\ref{eq:stimeW}) and classical Schauder estimates, we can conclude that $\sup_{v\in\MM_{2}}\|v\|_{\CC^{2}(\R^{2})^{2}}<+\infty$.}\end{Remark}

\begin{Remark}\label{R:unif} {\rm For every $v\in\MM_{2}$ we have $|v(x,y)-\a_{\pm}|\to 0$  as ${x\to\pm\infty}$ uniformly w.r.t. $y\in\R$. 
\noindent Indeed,
assume by contradiction that $v\in\MM_{2,q}$ and there 
exist $\delta>0$ and
a sequence $(x_{n},y_{n})\in\R^{2}$ with $x_{n}\to +\infty$ such 
that $|v(x_{n},y_{n})-\a_{+}|\geq 2\d$ for all $n\in\N$. 
By Remark \ref{R:trapped} there exists $\rho>0$ 
such that 
$
|v(x,y)-\a_{+}|\geq \d$ for all $(x,y)\in 
B_{\rho}(x_{n},y_{n})$, $n\in\N$. 
If along a subsequence we have $y_{n}\to 
y_{0}$,  then
$|v(x,y)-\a_{+}|\geq \d$ for all $(x,y)\in 
B_{\frac{\rho}{2}}(x_{n},y_{0})$ and $n$ large, which is not possible 
since $v(\cdot,y)\in \HH_1=z_{0}+H^{1}(\R)^{2}$ for almost every $y\in \R$. If $y_{n}\to +\infty$ (analogous to $y_{n}\to -\infty$) along a subsequence, we plainly obtain $
\limsup_{y\to +\infty}\|v(\cdot,y)-q\|_{L^{2}}>0
$
which contradicts $v\in \MM_{2,q}$.}\end{Remark}

\medskip

By definition,  we have that if $v\in\MM_{2,q}$ then
$|v(\cdot,y)-q|\to 0$ as $y\to+\infty$ with respect to the $L^{2}$ metric. We can in fact say more

\begin{Lemma}\label{L:unif}
If $v\in\MM_{2,q}$ for some  $q\in\MM_1^{*}$, then
$
\|v(\cdot,y)-q\|_{L^{\infty}}\to 0$
$\hbox{ as $y\to+\infty$}.
$
\end{Lemma}
{ \Proof
If $v\in\MM_{2,q}$, we know that $v$ solves ($P_2$) and by Remark \ref{R:trapped} that 
$\|v\|_{\CC^{2}}<+\infty$. 
If by contradiction $\limsup_{y\to+\infty}\|v(\cdot,y)-q\|_{L^{\infty}}\geq 2\rho_{0}>0$ then, there exists a sequence $(x_{n},y_{n})\in\R^{2}$ with $y_{n}\to +\infty$ such that $|v(x_{n},y_{n})-q(x_{n})|\geq\rho_{0}$. Then, since
$\|v-q\|_{\CC^{2}}<+\infty$, there exists $r_{0}>0$ such that
$|v(x,y_{n})-q(x)|\geq\rho_{0}/2$ whenever $|x-x_{n}|\leq r_{0}$ and $n\in\N$.
Hence we get the contradiction $\|v(\cdot,y_{n})-q\|_{L^{2}}\geq
r_{0}\rho_{0}^{2}/2$ for all $n\in\N$.\QED

\medskip

\noindent As last step in this section we use condition $(**)$ to obtain $L^{2}$ compactness of the minimizing sequences in $\tilde\HH_{2,q}$.\\ 
 
\begin{Lemma}\label{L:convL2}
Let $q\in\MM_1^{*}$ and $(v_{n})$ be a minimizing sequence in $\tilde\HH_{2,q}$ with $\|v_{n}\|_{L^{\infty}}\leq R$. Then 
$\dist_{L^{2}}(v_{n},\MM_{2,q})\to 0$ as $n\to +\infty$.
\end{Lemma}
\Proof 
Let $(v_{n})\subset\tilde\HH_{2,q}$ be such that $\varphi_{2}(v_{n})\to m_{2}$.  
By Lemma \ref{L:compattezza} and arguing as in the proof of Proposition \ref{P:sol2dim}, given any subsequence of $(v_{n})$, there exists $v\in\MM_{2,q}$ and a sub-subsequence, denoted again by $(v_{n})$, such that  $v_{n}- v\to 0$ weakly 
in $H^{1}(\R\times(-L,L))^{2}$ for any $L>0$. To prove the Lemma it will be sufficient to show that $\|v_{n}-v\|_{L^{2}}\to 0$.\\
Firstly note that, if $F$, $F_{1}$ and $F_{2}$ are   lower semicontinuous functionals such that $F=F_{1}+F_{2}$, if $v_{n}\to v$ and
$F(v_{n})\to F(v)$ then also $F_{1}(v_{n})\to F_{1}(v)$ and 
$F_{2}(v_{n})\to F_{2}(v)$. Iteratively applying this property in our case, since $\f(v_{n})\to\f(v)$,
we obtain that for any $L, M>0$, as $n\to+\infty$, we have
\begin{align}
\label{eq:semiintW}&\int_{-L}^{L}\int_{M}^{+\infty}W(v_{n})\, dx\, dy\to\int_{-L}^{L}\int_{M}^{+\infty}W(v)\, dx\, dy\\
\label{eq:semiint}
&\int_{L}^{+\infty}\varphi_1(v_{n}(\cdot,y))-m\ dy\to\int_{L}^{+\infty}\varphi_{1}( v(\cdot,y))-m\ dy
\end{align}
We use (\ref{eq:semiintW}) to deduce that 
 \begin{equation}\label{eq:convcomp}v_{n}- v\to 0
\hbox{  strongly
in }L^{2}(\R\times(-L,L))^{2}\hbox{  for any }L>0.\end{equation} 
Indeed, given $L>0$ and $\varepsilon>0$ we  fix $M>0$ such that $|v(x,y)-{\bf a}_{+}|<1$ for $x>M$ and $|y|<L$ and moreover  $\int_{-L}^{L}\int_{M}^{+\infty}W(v)\, dx<\varepsilon$. By (\ref{eq:semiintW}), there exists $\bar n\in\N$ such that $\int_{-L}^{L}\int_{M}^{+\infty}W(v_{n})\, dx\, dy<2\varepsilon$ for $n\geq\bar n$. Using (\ref{BGS}), we recover that there exists $C>0$ such that $\chi(\xi)\leq CW(\xi)$ for any $|\xi|\leq R$ so that, since $\|v_{n}\|_{L^{\infty}},\, \|v\|_{L^{\infty}}\leq R$, we deduce that for any $x>M$ and $|y|<L$ there results
$|v_{n}(x,y)-v(x,y)|^{2}\leq 2(\chi(v_{n}(x,y))^{2}+\chi(v(x,y))^{2})\leq 2C(W(v_{n}(x,y))+W(v(x,y)))$.
Then
$\|v_{n}-v\|^{2}_{L^{2}(\{|x|>M,\, |y|\leq L\})^{2}}\leq 4C \int_{-L}^{L}\int_{|x|>M}^{+\infty}W(v_{n})+W(v)\, dx\, dy\leq 12C\varepsilon$,
and since $\varepsilon$ is arbitrary and $v_{n}\to v$ in $L^{2}([-M,M]\times[-L,L])^{2}$, (\ref{eq:convcomp}) follows.\\
To conclude the proof  we show that
\begin{equation}\label{eq:claim}
\fa\,\e>0\, \exists\, L_{\e}>0,\, \bar n\in\N\hbox{ such that }
\|v_{n}-v\|_{L^{2}(\R\times\R\setminus {(-L_{\e}},L_{\e}))^{2}}\leq \e,\ \fa\, n\geq \bar n. 
\end{equation}
 It is not restrictive to assume that
$\f(v_{n})\leq m_{2}+\bar \lambda$ for any $n\in\N$, so that, by Lemma \ref{L:concentrazione2d} we have 
\begin{equation}\label{eq:minimoq-} 
\|v_{n}(\cdot,y)-q\|_{L^{2}}\leq d_{0}\hbox { for all }n\in\N\hbox{ and }y\geq\bar\ell
\end{equation} 
Letting $\varepsilon>0$, we fix $L_{\varepsilon}\geq\bar \ell$   such that $\int_{L_{\varepsilon}}^{+\infty}\F_{1}(v(\cdot,y))-m_{1}\ dy <\frac\varepsilon2$ and, by (\ref{eq:semiint}),  we fix also $n_{\varepsilon}\in\N$ such that 
\begin{equation}\label{eq:zerocomp}
\int_{L_{\varepsilon}}^{+\infty}\varphi_1(v_{n}(\cdot,y))-m_{1}\ dy<\varepsilon\hbox{ for all }n\geq n_{\varepsilon}.
\end{equation}
Then, letting $\nu^{*}>0$ as in
Lemma \ref{L:Fmenoc} and  denoting
$
\AA_{n}=\{ y>L_{\varepsilon}\,|\, \varphi_1(v_{n}(\cdot,y))-m_{1}> \nu^{*}\},$ 
 by (\ref{eq:zerocomp})
we have $\hbox{meas}(\AA_{n})\leq \frac{\varepsilon}{\nu^{*}}\hbox{ for any }n\geq n_{\varepsilon}$. Since $L_{\varepsilon}\geq\bar\ell$, by (\ref{eq:minimoq-}),
\begin{equation}\label{eq:primacomp}
\int_{\AA_{n}}\|v_{n}(\cdot,y)-q\|_{L^{2}}^{2}dy\leq \frac{\varepsilon d_{0}^{2}}{\nu^{*}}
\hbox{ for any }n\geq n_{\varepsilon},
\end{equation}
while,
since by  Lemma \ref{L:Fmenoc}   for $y\in(L_{\varepsilon},+\infty)\setminus\AA_{n}$ and $n\geq n_{\varepsilon}$ we have
\[
\|v_{n}(\cdot,y)-q\|_{L^{2}}^{2}=\dist_{L^{2}} (v_{n}(\cdot,y),\MM_1)^{2}\leq  \frac{4}{\omega^{*}}(\varphi_1(v_{n}(\cdot,y))-m).
\]
we recover
 \begin{equation}\label{eq:secondacomp}
\int_{(L_{\varepsilon},+\infty)\setminus\AA_{n}}\|v_{n}(\cdot,y)-q\|_{L^{2}}^{2}dy\leq \frac{4\varepsilon}{\omega^{*}}\hbox{ for any }n\geq n_{\varepsilon}.
\end{equation}
By (\ref{eq:primacomp}) and (\ref{eq:secondacomp}) for every $n\geq n_{\e}$ we obtain
$
\int_{(L_{\varepsilon},+\infty)}\|v_{n}(\cdot,y)-q\|_{L^{2}}^{2}dy\leq \varepsilon(\frac{4}{\omega^{*}}+\frac{d_{0}^{2}}{\nu^{*}})
$
and, by semicontinuity, the same estimate holds with $v_{n}$ replaced by $v$. Hence 
$
\int_{(L_{\varepsilon},+\infty)}\|v_{n}(\cdot,y)-v(\cdot,y)\|_{L^{2}}^{2}dy\leq 2\varepsilon(\frac{4}{\omega^{*}}+\frac{d_{0}^{2}}{\nu^{*}})\hbox{ for any }n\geq n_{\varepsilon}
$
and (\ref{eq:claim}) follows by the symmetry of $v$ and $v_{n}$.
\QED

\begin{Remark}\label{R:distM}
{\rm By (\ref{eq:d0bis}) and Lemma \ref{L:concentrazione2d} we have that $\dist_{L^{2}(\R\times(-L,L))^{2}}(\MM_{2,q},\MM_{2,p})\geq3d_{0}$ for all $p\not=q\in\MM_1^{*}$ and $L>\bar\ell+1$.\\ Indeed, if $v\in\MM_{2,q}$ and $w\in\MM_{2,p}$, by Lemma \ref{L:concentrazione2d} we have
$
\|p-q\|_{L^{2}}\le2d_{0}+\|v(\cdot,y)-w(\cdot,y)\|_{L^{2}}
$ for all $y\geq\bar\ell$,
and so, by (\ref{eq:d0bis}), $\|v(\cdot,y)-w(\cdot,y)\|_{L^{2}}\geq 3d_{0}$ for all $y\geq\bar\ell$. Then, for all $L>\bar\ell+1$ we conclude
$
\|v-w\|_{L^{2}(S_{L})^{2}}^{2}=\int_{-L}^{L}\|v(\cdot,y)-w(\cdot,y)\|_{L^{2}}^{2}dy\geq\int_{\bar\ell}^{\bar\ell+1}\|v(\cdot,y)-w(\cdot,y)\|_{L^{2}}^{2}dy\geq 9d_{0}^{2}.
$

}
\end{Remark}

\section{Two dimensional periodic solutions}

Here below, given $L>0$, we will study some variational properties of the minimal solutions to the problem 
\begin{equation*}\label{eq:solSL}
\begin{cases} -\Delta v(x,y)+\nabla W(v(x,y))=0,&
(x,y)\in S_{L},\cr 
v(x,y)=-v(-x,-y),& (x,y)\in S_{L},\cr
 \partial_{y} v(x,\pm L)={\bf 0},& x\in\R,\cr
 {\displaystyle\lim_{x\to\pm\infty}}|v(x,y)-\a_{\pm}|=0, 
	&\hbox{for every}y \in (-L,L),\cr\end{cases}\eqno(P_{L,2})
 \end{equation*} 
 where we have denoted $S_{L}=\R\times(-L,L)$.
Given $L>{0}$ we consider on the space 
$
\tilde\HH_{2,L}=\{v_{|S_{L}}\,|\, v\in\tilde\HH_{2}\}
$
the functional 
\begin{align*}
\varphi_{2,L}(v)\equiv\varphi_{2,(-L,L)}	&= \int_{-L}^{L}\tfrac12\|\partial_{y}v(\cdot,y)\|_{L^{2}}^{2}+\varphi_1(v(\cdot,y))-m_{1}\, dy
\end{align*}
which is positive on $\tilde\HH_{2,L}$ and  weakly lower semicontinuous w.r.t. the $H^{1}(S_{L})^{2}$ topology.
	We look for minima of $\varphi_{2,L}$ on $\tilde\HH_{2,L}$ and we set
	$$
m_{2,L}=\inf_{\tilde\HH_{2,L}}\varphi_{2,L}
\quad\hbox{and}\quad\MM_{2,L}=\{ v\in\tilde\HH_{2,L}\,|\, \varphi_{2,L}(v)= m_{2,L}\}
$$
noting that the map $L\mapsto m_{2,L}$ is not decreasing with $m_{2,L}\le  m_{2}$ for all $L>0$. 
\\

Using (\ref{eq:limitatezza}) and (\ref{eq:stimal2}), the argument in the proof of Lemma \ref{L:compattezza} applies to recover coerciveness properties of $\varphi_{2,L}$  analogous to the ones of $\f$. These properties allow us to apply the direct method of the Calculus of Variations to obtain

\begin{Proposition}\label{T:ML2}
For every $L>{0}$,  if $(v_{n})$ is a minimizing sequence for $\varphi_{2,L}$ in $\tilde\HH_{2,L}$, then there exists $\bar v\in\MM_{2,L}$ such that, up to a subsequence, $v_{n}\to \bar v$ weakly in $H^{1}(S_{L})^{2}$. In particular $\MM_{2,L}\not=\emptyset$ for every $L>0$.
\end{Proposition}
\begin{Remark}\label{R:solper}{\rm Note that, by ($W_{2}$), if $v\in\MM_{2,L}$ then $\|v\|_{L^{\infty}(S_{L})}\leq R$. Moreover, using the symmetry of $W$, it can be proved that every $v\in\MM_{2,L}$ is a classical solution to $(P_{L,2})$. From every $v\in\MM_{2,L}$ we then obtain a two dimensional periodic solution of ($P_{2}$) by  reflection with respect to $y=\pm L$ and then by periodic continuation in the $y$-direction. We will still denote with $v$ the obtained solution on $\R^{2}$ noting that it is $y$-periodic of period $4L$ and satisfies $\partial_{y}v(x,\pm L)={\bf 0}$ for all $x\in\R$. Again by ($W_{2}$), since $\|v\|_{L^{\infty}(S_{L})}\leq R$, Schauder estimates give the existence of a constant $C>0$ such that
\begin{equation}\label{eq:stimeC2}
\|v\|_{\CC^{2}(\R^{2})^{2}}\leq C\hbox{ for all }L>0\hbox{ and }v\in\MM_{2,L}.
\end{equation}}\end{Remark}
\begin{Lemma}\label{L:contm2L} The map $L\mapsto m_{2,L}$ is continuous on $\R^{+}$. Moreover if $L_{n}\to L_{0}>0$ and $v_{n}\in\MM_{2,L_{n}}$ for $n\in\N$, then there exists $\bar v\in\MM_{2,L_{0}}$ such that, up to a subsequence, $v_{n}- \bar v\to 0$ weakly in $H^{1}(S_{L_{0}})^{2}$ and strongly in $L^{2}(S_{L_{0}})^{2}$.\end{Lemma} 
\Proof Given $L_{0}>0$, $L_{n}\to L_{0}^{+}$ and $v_{0}\in\MM_{2,L_{0}}$ we have $m_{2,L_{n}}\leq\varphi_{2,L_{n}}(v_{0})\to m_{2,L_{0}}$ as $n\to+\infty$, so that, by monotonicity of $m_{2,L}$, we derive  $m_{2,L_{n}}\to m_{2,L_{0}}$. Hence, if $v_{n}\in\MM_{2,L_{n}}$ we obtain $\varphi_{2,L_{0}}(v_{n})\to m_{2,L_{0}}$ and, by Proposition \ref{T:ML2} there exists $\bar v\in\MM_{2,L_{0}}$ such that, up to a subsequence, $v_{n}\to \bar v$ weakly in $H^{1}(S_{L_{0}})^{2}$.\\
Consider now a sequence $L_{n}\to L_{0}^{-}$ and let $v_{n}\in\MM_{2,L_{n}}$. Since  $\varphi_{2,L_{0}}(v_{n})\leq C$ and $\|v_{n}\|_{L^{\infty}}\leq R$, by  (\ref{eq:limitatezza}) and (\ref{eq:stimal2}) we recover that $(v_{n}-z_{0})$ is bounded in
$H^{1}(S_{L_{0}})^{2}$ and, up to a subsequence,
it converges, weakly in $H^{1}(S_{L_{0}})^{2}$ and almost everywhere, to a function $\bar v-z_{0}$ with $\bar v\in \tilde\HH_{2,L_{0}}$.
By semicontinuity, for any $L\in (0, L_{0})$, we have $\varphi_{2,L}(\bar v)\leq\liminf_{n\to+\infty}\varphi_{2,L}(v_{n})\leq\lim_{n\to+\infty}m_{2,L_{n}}$
and we derive
$m_{2,L_{0}}\leq\varphi_{2,L_{0}}(\bar v)=\sup_{L<L_{0}}\varphi_{2,L}(v)\leq \lim_{n\to+\infty}m_{2,L_{n}}\leq m_{2,L_{0}}$, from which $m_{2,L_{n}}\to m_{2,L_{0}}$ and $\bar v\in\MM_{2,L_{0}}$.\\
The strong convergence in $L^{2}(S_{L_{0}})^{2}$ follows by applying the semicontinuity argument used to derive (\ref{eq:convcomp}) in the proof of Lemma \ref{L:convL2}.
\QED

\noindent We now study the behaviour of $\MM_{2,L}$ and $m_{2,L}$ for $L\to+\infty$.
Given $q\in\MM_1$, for $L>\bar\ell$, we  define    
    $$
    \tilde\HH_{2,L,q}=
    \{v\in\tilde\HH_{2,L}\,|\, 
    \|v(\cdot,y)-q\|_{L^{2}}\leq d_{0},\,\forall
		y\in[\bar\ell,L]\}
    $$
    and we note that by Lemma \ref{L:concentrazione2dL}, there results $\MM_{2,L}\subset\cup_{q\in\MM_1^{*}}\tilde\HH_{2,L,q}$ for every $L\geq \bar\ell$. We denote $\MM_{2,L,q}=\MM_{2,L}\cap\tilde\HH_{2,L,q}$.
    
 \begin{Lemma}\label{L:C2Ltoinfty}
 There results $m_{2,L}\to m_{2}$ as $L\to+\infty$.
 \end{Lemma}
 \Proof
 Let $L_{n}\to+\infty$ and $v_{n}\in\MM_{2,L_{n}}$. By the above considerations, there exists $q\in \MM_{1}^{*}$ such that along a subsequence, denoted again $v_{n}$, we have $v_{n}\in\MM_{2,L_{n},q}$.
 Since $(\F_{2,(\bar\ell,L_{n})}(v_{n}))_{n\in\N}$ is bounded and since $L_{n}\to+\infty$ we recover that there exists $y_{n}\in
 (\bar\ell,L_{n})$ such that $\F_{1}(v_{n}(\cdot,y_{n}))\to m_{1}$. By Lemma \ref{L:compK}, since $v_{n}\in\tilde\HH_{2,L_{n},q}$ and so $\|v_{n}(\cdot,y_{n})-q\|_{L^{2}}\leq r_{0}$, we obtain that, 
 $\|v_{n}(\cdot,y_{n})-q\|_{H^{1}}\to 0$. We define 
 $$
\bar v_{n}(x,y)=\begin{cases}   q(x)&  y\geq y_{n}+1, \cr           
	(y_{n}+1-y)v_{n}(x,y_{n})+(y-y_{n})q(x) & y_{n}\le y\le y_{n}+1\cr
	v_{n}(x,y)& -y_{n}\le y\le y_{n}\cr
	(y_{n}+1+y)v_{n}(x,-y_{n})+(-y-y_{n})\bar q(x) &-y_{n}\geq y\geq -y_{n}-1\cr
	\bar q(x)&  y\le -y_{n}-1,
	\end{cases}
$$
observing  that $\bar v_{n}\in\tilde\HH_{2,q}\cap\tilde\HH_{2,q,L_{n}}$  so that $\F_{2,L_{n}}(\bar v_{n})\geq m_{2,L_{n}}$.
We have 
$
\|\bar v_{n}(\cdot,y)-q\|_{H^{1}}\leq \|\bar v_{n}(\cdot,y_{n})-q\|_{H^{1}}$ and  
$\|\partial_{y}\bar v_{n}(\cdot,y)\|_{L^{2}}=\|v_{n}(\cdot,y_{n})-q\|_{L^{2}}
$ for any $y\in ( y_{n},y_{n}+1)$
and since $\|v_{n}(\cdot,y_{n})-q\|_{H^{1}}\to 0$, we derive that 
 $\F_{1}(\bar v_{n}(\cdot,y))\to m_{1}$
uniformly w.r.t. $y\in ( y_{n},y_{n}+1)$} and then
 $\F_{2,( y_{n}, y_{n}+1)}(\bar v_{n})\to 0$. Hence $\F_{2,(y_{n},L_{n})}(\bar v_{n})\to 0$ too and so
 $\F_{2,(-y_{n},y_{n})}(\bar v_{n})=\F_{2,L_{n}}( v_{n})+o(1)=m_{2,L_{n}}+o(1)$. Therefore
 $m_{2}\leq\F_{2}(\bar v_{n})=\F_{2,(-y_{n},y_{n})}(v_{n})+2\F_{2,(y_{n}, y_{n}+1)}(\bar v_{n})=m_{2,L_{n}}+o(1)\leq m_{2}+o(1),$ and the Lemma follows.
\QED
\begin{Remark}\label{R:tend0}{\rm The argument in the proof of Lemma \ref{L:C2Ltoinfty} can be used to show that if $q\in\MM_{1}^{*}$, $L_{n}\to+\infty$, $y_{n}\in (\bar\ell,L_{n})$ and $v_{n}\in\tilde\HH_{2,L_{n},q}$ are such that $\F_{2,L_{n}}(v_{n})-m_{2,L_{n}}\to 0$ and $\varphi_{1}(v_{n}(\cdot,y_{n}))\to m_{1}$
then $\F_{2,(y_{n},L_{n})}(v_{n})\to0$ as $n\to +\infty$.\\ 
As a consequence, note that if $q\in\MM_{1}^{*}$, $v_{n}\in\tilde\HH_{2,L_{n},q}$ $\F_{2,L_{n}}(v_{n})-m_{2,L_{n}}\to 0$ and $\ell_{n}\in(\bar\ell,L_{n})$ is such that $\ell_{n}\to+\infty$, then, since
there exists $y_{n}\in
 (\bar\ell,\ell_{n})$ such that $\F_{1}(v_{n}(\cdot,y_{n}))\to m_{1}$, we obtain $\F_{2,(\ell_{n},L_{n})}(v_{n})\to 0$.}\end{Remark}

\begin{Lemma}\label{L:compattezzaperiodiche}
 Let $q\in\MM_1^{*}$, $L_{n}\to+\infty$ and $v_{n}\in\tilde\HH_{2,L_{n},q}$ be such that 
 $\F_{2,L_{n}}(v_{n})-m_{2,L_{n}}\to 0$ as $n\to\infty$. Then 
 $
 \dist_{L^{2}(S_{L_{n}})^{2}}(v_{n},\MM_{2,q})\to 0
 $  as $n\to +\infty$.
  \end{Lemma} 

 \Proof
Fix any $\varepsilon>0$. By Remark \ref{R:tend0} we have that $\F_{2,(L_{n}-\varepsilon,L_{n})}(v_{n})\to 0$ as $n\to+\infty$ and so there exists  $y_{n}\in
 (L_{n}-\varepsilon,L_{n})$ such that $\F_{1}(v_{n}(\cdot,y_{n}))\to m_{1}$. We argue as
 in the proof of Lemma \ref{L:C2Ltoinfty} defining the function $\bar v_{n}$ relatively to the sequence $y_{n}$, obtaining that
$
\F_{2}(\bar v_{n})\to m_{2}\hbox{ as }n\to\infty.
$
By Lemma \ref{L:convL2} we conclude that
$
\dist_{L^{2}}(\bar v_{n},\MM_{2,q} )\to 0
$,
and  so, since $\bar v_{n}(x,y)=
v_{n}(x,y)$ on $S_{y_{n}}$, we derive
$
\dist_{L^{2}(S_{y_{n}})^{2}}(v_{n},\MM_{2,q})\to 0
$. By Lemma \ref{L:concentrazione2dL} we know that if $n$ is sufficiently large, $\|v_{n}(\cdot,y)-q\|_{L^{2}}\leq d_{0}$ for $y\in(\bar\ell,L_{n})$, 
and we conclude that $\dist_{L^{2}(S_{L_{n}})^{2}}(v_{n},\MM_{2,q})=o(1)+2\varepsilon d_{0}$. The Lemma then follows by the arbitrariness of $\varepsilon$.\QED

As a direct consequence of Lemmas \ref{L:C2Ltoinfty} and \ref{L:compattezzaperiodiche} we obtain
\begin{Corollary}\label{C:mu0L2}
 There exists $\bar L\geq\bar\ell$ such that for all $\delta>0$ there exists $\mu_{\delta}>0$ such that 
 if $L\geq  \bar L$, $v\in\tilde \HH_{2,L}$ and $\min_{q\in\MM_1^{*}}\dist_{L^{2}(S_{L})^{2}}(v,\MM_{2,q})>\delta$ then $\F_{2,L}(v)\geq m_{2,L}+\mu_{\delta}$.
 \end{Corollary}
 
Next Lemma describes the uniformity of the behaviour of the functions in $\MM_{2}\cup \cup_{L>\bar\ell}\MM_{2,L}$ when $|x|$ is large.

\begin{Lemma}\label{L:concLinftyxlarge} For all $\varepsilon>0$ there exists $T_{\varepsilon}>0$ such that if 
$v\in\MM_{2}\cup \cup_{L\geq\bar\ell}\MM_{2,L}$, then
$|v(x,y)-\a_{+}|<\e$
for all $y\in \R$ and $x>T_{\varepsilon}$.
\end{Lemma}
\Proof The argument in Remark \ref{R:unif} can be used to prove that for any $L>0$ and $v\in\MM_{2,L}$ we have
$|v(x,y)-\a_{\pm}|\to 0$  as ${x\to\pm\infty}$, uniformly for $y\in [-L,L]$ (and so for $y\in\R$ by periodicity).
To obtain the uniform result, we argue by contradiction, assuming that there exists $\delta>0$,  $(L_{n})\subset [\bar\ell,+\infty)$, $v_{n}\in\MM_{2,L_{n}}$ and  $(x_{n},y_{n})\in S_{L_{n}}$ such that $x_{n}\to+\infty$ and  $|v_{n}(x_{n},y_{n})-\a_{+}|\geq 2\d$ for all $n\in\N$.
By Remarks \ref{R:trapped} and \ref{R:solper} there exists $\rho>0$ 
such that 
$
|v_{n}(x,y)-\a_{+}|\geq \d$ for all $(x,y)\in 
B_{\rho}(x_{n},y_{n})$, $n\in\N$. If $L_{n}$ accumulates some $L_{0}\geq\bar\ell$, by Lemma \ref{L:contm2L} there exists
$\bar v\in\MM_{2,L_{0}}$ such that, along a subsequence, $v_{n}- \bar v\to 0$ strongly in $L^{2}(S_{L_{0}})^{2}$.
Then $\|\bar v-a_{+}\|_{L^{2}(B_{\rho}(x_{n},y_{n}))^{2}}\geq  \|v_{n}-a_{+}\|_{L^{2}(B_{\rho}(x_{n},y_{n}))^{2}}+o(1)\geq \pi\rho^{2}\d+o(1)$ in contradiction with the fact that $|v(x,y)-\a_{+}|\to 0$  as ${x\to\pm\infty}$ uniformly for $y\in (-L_{0},L_{0})$.
The case $L_{n}\to +\infty$ and the case $v\in\MM_{2}$, leads to an analogous contradiction using Lemmas \ref{L:convL2} and \ref{L:compattezzaperiodiche}.
\QED

Moreover we have
\begin{Lemma}\label{L:unif2L} There exists $\bar T>0$ such that for any $v\in \MM_{2}\cup \cup_{L\geq\bar\ell}\MM_{2,L}$ we have
$|v(x,y)-\a_{+}|\leq\bar \d \hbox{e}^{-\sqrt{\underline {w}/2}(x-\bar T)}$ for $x\geq \bar T$ and $y\in\R$.\end{Lemma}
\Proof
Let $\bar\delta$ be defined by Remark \ref{R:costanti} and let $\bar T\geq\max\{ T_{\bar\delta},T_{0}\}$, where $T_{\bar\d}$ is given by Lemma \ref{L:concLinftyxlarge} for $\varepsilon=\bar\d$ and $T_{0}$ by Remark  \ref{R:boundLinfty+conc}. By Lemma  \ref{L:concLinftyxlarge} we have that $|v(x,y)-\a_{+}|< \bar\delta$ for any
$x>\bar T$ and $v\in \MM_{2}\cup \cup_{L\geq\bar\ell}\MM_{2,L}$. Defining $\psi(x,y)=|v(x,y)-\a_{+}|^{2}$, by (\ref{eq:iper}) we have 
$\Delta \psi\geq 2\underline{w}\psi$  for $x>\bar T$. Considering the function $\xi(x)=
\bar\delta^{2}\hbox{e}^{-\sqrt{2\underline{w}}(x-\bar T)}$ we have that
$\Delta (\psi-\xi)\geq 2\underline{w}(\psi-\xi)$. \\
We have  that $\psi-\xi<0$. Indeed, if $v\in\MM_{2,L}$ for a $L\geq\bar\ell$, then
the function $\psi-\xi$ is periodic in the variable $y$, $(\psi-\xi)(\bar T,y)<0$ and $\lim_{x\to+\infty}(\psi-\xi)(x,y)=0$
 for $y\in\R$. If $\psi-\xi$ assumes a positive value on $\{(x,y)\subset\R^{2}\,/\, x\geq \bar T\}$ we recover that it has a positive maximum 
 on $\{(x,y)\subset\R^{2}\,/\, x\geq \bar T\}$ which is not possible since $\Delta (\psi-\xi)\geq 2\underline{w}(\psi-\xi)$.\\
 If $v\in\MM_{2}$, by Lemma \ref{L:unif} there exists $q\in\MM_{1}^{*}$ such that $\|v(\cdot,y)-q\|_{L^{\infty}}\to 0$ as $y\to +\infty$. By Lemma \ref{L:expSOL1dim}, since
 $\bar T\geq T_{0}$ we now that $|q(x)-\a_{+}|^{2}\leq \bar\delta^{2}\hbox{e}^{-\sqrt{2\underline{w}}(x-\bar T)}$ and so by Lemma
 \ref{L:unif} and by symmetry we have that $\limsup_{y\to\pm\infty}\psi(x,y)-\xi(x)\leq 0$ for any $x\geq \bar T$. Since $\psi(\bar T,y)-\xi(\bar T)\leq 0$ and
  $\psi(x,y)-\xi(x)\to  0$ as $x\to +\infty$, we conclude $\psi-\xi\leq 0$ on $\{(x,y)\subset\R^{2}\,/\, x\geq \bar T\}$ since
  otherwise $\psi-\xi$ has a positive maximum 
 on $\{(x,y)\subset\R^{2}\,/\, x\geq \bar T\}$ which is again not possible.\QED

We can also uniformly characterize the asymptotic behaviour for $y$ large of the functions in $\MM_{2,L}$. First we give a ''$L^{2}$-concentration'' result.

\begin{Lemma}\label{L:concLn}
Let $q\in\MM_1^{*}$, $L_{n}\to+\infty$ and $v_{n}\in\tilde\HH_{2,L_{n},q}$ be such that 
 $\F_{2,L_{n}}(v_{n})-m_{2,L_{n}}\to 0$ as $n\to\infty$. Then for any $\varepsilon>0$ there exists
 $\ell_{\varepsilon}\geq\bar\ell$ and $n_{\varepsilon}\in\N$ such that $\sup_{y\in (\ell_{\varepsilon},L_{n})}\|v_{n}(\cdot,y)-q\|_{L^{2}}<\varepsilon$ for any $n\geq n_{\varepsilon}$.\end{Lemma}
\Proof
Assume by contradiction that there exist $\delta>0$ and a sequence $\ell_{n}\in (\bar\ell,L_{n})$ such that $\ell_{n}\to +\infty$ and for which  $\|v_{n}(\cdot,\ell_{n})-q\|_{L^{2}}\geq 2\delta$. Since $\ell_{n}\to +\infty$ there exists a sequence $y_{n}\in (\ell_{n}/2,\ell_{n})$ such that $\varphi_{1}(v_{n}(\cdot,y_{n}))-m_{1}\to0$ and so, by Lemma \ref{L:compK}, we have that $\|v_{n}(\cdot,y_{n})-q\|_{H^{1}}\to 0$. 
 Then, by (\ref{eq:contd}),
there exists $(\sigma_{n},\tau_{n})\subset(y_{n},\ell_{n})$ such that $\|v_{n}(\cdot,\tau_{n})-
v_{n}(\cdot,\sigma_{n})\|_{L^{2}}=\delta$ and $\dist_{L^{2}}(v_{n}(\cdot,y),\MM_1)_{L^{2}(\R)}\geq \delta$ for  $y\in
(\sigma_{n},\tau_{n})$. Hence by (\ref{eq:h0}) there exists $\nu>0$ such that
$\F_{1}(v_{n}(\cdot,y))\geq m_{1}+\nu$ for almost every $y\in
(\sigma_{n},\tau_{n})$ and by (\ref{eq:stime2dim}) we conclude
$\F_{2,(\sigma_{n},\tau_{n})}(v_{n})\geq\sqrt{2\delta}\nu>0$ for any $n\in\N$, in contradiction with Remark \ref{R:tend0}.
\QED 

Specializing Lemma \ref{L:concLn} to functions in $\MM_{2,L}$ we obtain also $L^{\infty}$ estimates.

\begin{Lemma}\label{L:concLinfty} For all $\varepsilon>0$ there exists $L_{\varepsilon}\geq\bar\ell$ such that if $L\geq L_{\varepsilon}$ and $v\in\MM_{2,L,q}$ for some $q\in\MM_1^{*}$, then
$
\|v(\cdot,y)-q\|_{L^{\infty}}<\e$
for all $y\in (L_{\varepsilon},L].
$
\end{Lemma}
\Proof
 Assume by contradiction  that there exist $\eta_{0}\in (0,\min\{\frac12;d_{0}\})$, $q\in\MM_{1}^{*}$,
two  sequences $0<\bar\ell<y_{n}\leq L_{n}$ such that $y_{n}\to+\infty$ 
and a sequence $v_{n}\subset\MM_{2,L_{n},q}$ such that
$\|v_{n}(\cdot,y_{n})-q\|_{L^{\infty}}>2\eta_{0}$  for $n\in\N$ . By (\ref{eq:stimeC2}) we have $\sup_{n\in\N}\|v_{n}-q\|_{\CC^{2}(S_{L_{n}})^{2}}<+\infty$ and so there exists 
$\delta\in (0,\bar\ell)$ such that $\|v_{n}(\cdot,y)-q\|_{L^{\infty}}>\eta_{0}$ for all $y\in(y_{n}-\delta,y_{n}).
$
Again by (\ref{eq:stimeC2}) we deduce that there esists $\e_{0}>0$ such that
$
\|v_{n}(\cdot,y)-q\|_{L^{2}}>\e_{0}$ for all $y\in(y_{n}-\delta,y_{n}),
$ in contradiction with Lemma \ref{L:concLn}.\QED

\noindent Assumption $(**)$ and Lemma \ref{L:Fmenoc} are now used to obtain exponential estimates. 

\begin{Lemma}\label{L:stimaexp2dim1}
There exist $\tilde L\geq \bar\ell$, $C>0$ and $\bar\omega>0$ such that if $L\geq\tilde L$ and $v\in\MM_{2,L,q}$ for some $q\in\MM_1^{*}$ then 
$\|v(\cdot,y)-q\|_{H^{1}}\leq  C\hbox{e}^{-\bar\omega y}$
for all 
$y\in[\tilde L,L]$
\end{Lemma}
\Proof  Given $v\in \MM_{2,L,q}$, consider the function 
$\phi(y)=\|v(\cdot,y)-q\|_{L^{2}}^{2}$ for $y> 0$. We have that
$\phi\in C^{2}((0,+\infty))$ and 
for every $y>0$ 
\begin{equation}\label{eq:sstima}
\ddot\phi(y)=2\|\partial_{y}v(\cdot,y)\|^{2}_{L^{2}}+
2\langle \partial_{yy}u(\cdot,y),v(\cdot,y)-q,\rangle_{L^{2}}\geq2 \langle \partial_{yy}v(\cdot,y),v(\cdot,y)-q\rangle_{L^{2}}\end{equation}
By Lemma \ref{L:concLinfty},  we can fix  $\tilde L\geq \bar\ell$ such that for  $y\in [\tilde L,L]$ and  $x\in\R$ we have
$
\nabla W(v(x,y))-\nabla W(q(x))= \nabla^{2}W(q(x))(v(x,y)-q(x))+r(x,y)
$
with $|r(x,y)|\leq \tfrac {\omega^{*}}{2}|
v(x,y)-q(x)|$ ($\omega^{*}$ defined in $(**)$). Since $v$ and $q$ solve (\ref{eq:eq0}) we have that
$\langle \partial_{yy}v(\cdot,y), v(\cdot,y)-q\rangle_{L^{2}}=\F''_{1}(q)(v(\cdot,y)-q)\cdot(v(\cdot,y)-q)+
\langle \nabla W(v(\cdot,y))-\nabla W(q)-\nabla^{2}W(q)(v(\cdot,y)-q), v(\cdot,y)-q\rangle_{L^{2}}\geq
\F''_{1}(q)(v(\cdot,y)-q)\cdot(v(\cdot,y)-q)-\tfrac {\omega^{*}}{2}\|v(\cdot,y)-q\|_{L^{2}}^{2}$ and so,
by $(**)$, 
$\langle \partial_{yy}v(\cdot,y),v(\cdot,y)-q\rangle\geq\tfrac {\omega^{*}}{2}\|v(\cdot,y)-q\|_{L^{2}}^{2}$
for every $y\in[\tilde L,L]$.
By (\ref{eq:sstima}) we obtain
\begin{equation}\label{eq:finalestime2dim}
\ddot\phi(y)\geq {\omega^{*}}\phi(y)
\hbox{ for all $y\in[\tilde L,L]$}.
\end{equation}
Since $v$ is symmetric with respect to $y=L$, the inequality (\ref{eq:finalestime2dim}) holds true also for all $y\in
(\tilde L,2L-\tilde L)$. Moreover, since $\tilde L\geq\bar\ell$ and $v\in\MM_{2,L,q}$,  there results $\phi(\tilde L)=\phi(2L-\tilde L)\le d_{0}^{2}$. Defining $\psi(y)=d_{0}^{2}\frac{\cosh(\sqrt{\omega^{*}}(L-y))}{\cosh(\sqrt{\omega^{*}}(L-\tilde L))}$, by (\ref{eq:finalestime2dim}) we have $(\ddot\phi-\ddot\psi)(y)\geq {\omega^{*}}(\phi-\psi)(y)$ on
$(\tilde L,2L-\tilde L)$ and $\phi-\psi\leq 0$ on $\partial (\tilde L,2L-\tilde L)$. Hence $\phi (y)\leq\psi(y)$
for $y\in [\tilde L,2L-\tilde L]$ and we conclude that
\begin{equation}\label{eq:stimaexp2dim1}
\|v(\cdot,y)-q\|_{L^{2}}^{2}\leq 
d_{0}^{2}\hbox{e}^{\sqrt{\omega^{*}}(\tilde L-y)}\quad\forall y\in[\tilde L,L].
\end{equation}
To conclude the proof note that, thanks to Lemmas \ref{L:expSOL1dim} and \ref{L:unif2L}, there exists $C_{1}>0$, independent of $L$, $v$ and $q$, such that $|v(x,y)-q(x)|\leq
C_{1}e^{-\sqrt{\underline{w}/2}\,x}$ for all $x>0$ and $y\in\R$. Then, setting $\bar\omega=\frac 12\min\{\sqrt{\omega^{*}}, \sqrt{\underline{w}/2}\}$, since for $y>0$ we have
$\|v(\cdot,y)-q\|_{L^{1}}=\int_{-y}^{y}|v(x,y)-q(x)|\, dx+2\int_{x\geq y}|v(x,y)-q(x)|\, dx\leq \sqrt{2y}\|v(\cdot,y)-q\|_{L^{2}}+4C\frac{\bar\d}{\sqrt{\underline{w}}}e^{-\sqrt{\underline{w}/2}y}$, by (\ref{eq:stimaexp2dim1}) we recover the existence of a constant $C_{2}>0$ such that $\|v(\cdot,y)-q\|_{L^{1}}\leq C_{2}e^{-\bar\omega y}$ for any $L\geq \tilde L$, $y\in[\tilde L,L]$.
 The Lemma then follows since
$
\|\partial_{x}(v(\cdot,y)-q)\|^{2}_{L^{2}}=-\int_{\R}(v(x,y)-q(x))\cdot\partial_{xx}(v(x,y)-q(x))\, dx$ and
since, by (\ref{eq:stimeC2}),  $|\partial_{xx}v(x,y)-\ddot q(x)|\le C$ on $\R^{2}$.
\QED

\noindent The last result of the section concerns the behaviour of  $m_{2,L}$ as $L\to +\infty$.

\begin{Lemma}\label{L:asimpotic-cL2}
There exists a constant $C>0$ such that for all $L\geq \tilde L$ there results
$m_{2}-m_{2,L}\leq C\hbox{e}^{-\frac{\bar\omega}2 L}$.
\end{Lemma}
\Proof Letting $v\in\MM_{2,L,q}$,
we define
$$
\tilde v(x,y)=\begin{cases}   q(x)&  y\geq L+1, \cr           
	(L+1-y)v(x,L)+(y-L)q(x) & L\le y\le L+1\cr
	v(x,y),& -L\le y\le L\cr
	(L+1+y)v(x,-L)+(-y-L)\bar q(x) &-L\geq y\geq -L-1\cr
	\bar q(x)&  y\le -L-1,
	\end{cases}
$$
Then $\tilde v\in\tilde\HH_{2,q}$ and $m_{2}\leq\F_{2}(\tilde v)= m_{2,L}+2\F_{2,(L,L+1)}(\tilde v)$. Hence
\begin{equation}\label{eq:livelli2} 
m_{2}-m_{2,L,q}\leq 2\F_{2,(L,L+1)}(\tilde v)=2\int_{L}^{L+1}\tfrac12\|\partial_{y}\tilde v(\cdot,y)\|^{2}_{L^{2}}+\varphi_1(\tilde v(\cdot,y))-m_{1}\,dy.
\end{equation}
We now  note that
by Lemma \ref{L:stimaexp2dim1} there exists a constant $C_{1}>0$ (independet of $L$) such that for every $y\in (L,L+1)$ there results
\begin{equation}\label{eq:livelli21} 
\|\partial_{y}\tilde v(\cdot,y)\|^{2}_{L^{2}}=\| v(\cdot,L)-q\|_{L^{2}}^{2}\leq C_{1}\hbox{e}^{-{2\bar\omega}L}.
\end{equation}
To give an estimate of $\F_{1}(\tilde v(\cdot,y))-m_{1}= \F_{1}(\tilde v(\cdot,y))-\F_{1}(q)$, we note first that  $
||\partial_{x} \tilde v(x,y)|^{2}-|\dot q(x)|^{2}|
\le |\partial_{x}v(x,L)-\dot q(x)|^{2}+2|\partial_{x}v(x,L)-\dot q(x)||\dot q(x)|$ for all $y\in (L,L+1)$. Applying Lemma \ref{L:stimaexp2dim1}
we then have that, for a certain constant $C_{2}>0$ (independent of $L$) and $y\in(L,L+1)$, there results
\[|\|\partial_{x}\tilde v(\cdot,y)\|^{2}_{L^{2}}-\|\dot q\|^{2}_{L^{2}}|
\le
 \|\partial_{x}v(\cdot,L)-\dot q\|_{L^{2}}^{2}+2\|\partial_{x}v(\cdot,L)-\dot q\|_{L^{2}}\|\dot q\|_{L^{2}}\le C_{2}e^{-\bar\omega L}\]%\\
Secondly observe that
 $|W(\tilde v(x,y))-W(q(x))|
\le \tilde w|v(x,L)-q(x)|$ for all $y\in(L,L+1)$ (where $\tilde w=\sup_{|\xi|<R}|\nabla W(\xi)|$)
and so, since we already know (as in the proof of the preceding Lemma)  that
$\|v(\cdot,L)-q\|_{L^{1}}\leq C_{3}e^{-\bar\omega L}$ for a certain constant $C_{3}>0$,
we conclude that $|\int_{\R}W(\tilde v(x,y))-W(q(x))dx|\leq C_{4}e^{-\bar\omega L}$ for any $y\in [L,L+1]$.
The above estimates allow us to conclude that there exists a constant $C_{5}>0$ such that  
$\varphi_1(\tilde v(\cdot,y))-m_{1}
\le C_{5}e^{-{\bar\omega}L}$ for every $y\in(L,L+1)$ and
the Lemma follows by
(\ref{eq:livelli21}) and (\ref{eq:livelli2}).\QED

\section{Three dimensional solutions}
The proof of Theorem \ref{T:main}, discussed in this last section, is based on a variational renormalized procedure developed for the scalar case in \cite{[AlCMsaddle]}, \cite{[AlMtran]} and  \cite{[AlMsadd3]}.\medskip 

Fixed $\theta\in(0,\frac\pi4]$ and denoted $\bar z=z\tan\theta$ and
$
\PP_{\theta}=\{(x,y,z)\in\R^{3}\,|\, (x,y)\in S_{\bar z},\, z\geq 0\}
$
we consider the space
$$
\HH_{3,\theta}=\{ u\in H^{1}_{loc}(\PP_{\theta})^{2}\,|\, u(\cdot,\cdot,z)\in\tilde\HH_{2,\bar z}\hbox{ for almost every } z\in(0,+\infty)\}
$$
on which we will look for minima of the functional
\begin{align*}
\varphi_{3}(u)&=\int_{0}^{+\infty}\int_{-\bar z}^{\bar z}\int_{\R}\tfrac12|\nabla u(x,y,z)|^{2}+W(u(x,y,z))\,dx-m_{1}\,dy-m_{2,\bar z}\, dz\\
&=\int_{0}^{+\infty}\tfrac12\|\partial_{z}u(\cdot,\cdot, z)\|_{L^{2}(S_{\bar z})^{2}}^{2}+\varphi_{2,\bar z}(u(\cdot, \cdot,z))-m_{2,\bar z}\,dz.
\end{align*}
Given an interval $I\subset\R_{+}$  we will also consider on $\HH_{3,\theta}$ the functional
$$
\varphi_{3,I}(u)=\int_{I}\tfrac12\|\partial_{z}u(\cdot,\cdot,z)\|^{2}_{L^{2}(S_{\bar z})^{2}}+
\F_{2,\bar z}u(\cdot,\cdot,z))-m_{2,\bar z}\, dz.
$$
We will denote
$$m_{3,\theta}=\inf_{u\in\HH_{3,\theta}}\F_{3}(u)\quad \hbox{and}\quad \MM_{3,\theta}=\{u\in\HH_{3,\theta}\,|\, \F_{3}(u)=m_{3,\theta}\}.
$$

\begin{Remark}\label{R:note3}{\rm If $u\in\HH_{3,\theta}$ then $u(\cdot,\cdot,z)\in\tilde\HH_{2,\bar z}$ for almost every $z>0$ and so $\F_{2,\bar z}(u(\cdot,\cdot,z))\geq m_{2,\bar z}$ for almost every $z>0$. Hence $\F_{3,I}$ is well defined on $\HH_{3,\theta}$  and non negative for every interval $I\subseteq \R^{+}$. It is moreover standard to show that
$\F_{3,I}$ is lower semicontinuous with respect to the weak topology of $H^{1}_{loc}(\PP_{\theta})^{2}$.}
\end{Remark}
\begin{Remark}\label{L:min3}
{\rm We have $m_{3}<+\infty$. Indeed, if $v\in\MM_{2,q}$ for a $q\in\MM_{1}^{*}$, the function $u(x,y,z)=v_{|\PP_{\theta}}(x,y)\in\HH_{3,\theta}$ and, by  Lemma \ref{L:asimpotic-cL2}, 
$
\F_{3}(u)=\int_{0}^{+\infty}\varphi_{2,\bar z}(v)-m_{2,\bar z}\,dz\leq \int_{0}^{+\infty}m_{2}- m_{2,\bar z}\, dz\leq C\int_{0}^{+\infty}\hbox{e}^{-\frac{\bar\omega}2 \tan\theta z}\, dz<+\infty.
$ }
\end{Remark}
\begin{Remark}\label{R:l3estimate}{\rm Denoted $T_{r}=\PP_{\theta}\cap\{z<r\}$ for  $r>0$, we have that if $u\in\HH_{3,\theta}$ then
\begin{align}
\nonumber\norm{\nabla u}^{2}_{L^{2}(T_{r})}&\leq 2\varphi_{3}(u)+2\int_{0}^{r}m_{2,\bar z}\, dz+2\int_{0}^{r}\int_{-\bar z}^{\bar z}m_{1}\,dy\,dz\\\label{eq:stima3a}
&\leq 2(\varphi_{3}(u)+rm_{2}+\tan(\theta)r^{2}m_{1})
\end{align}
Moreover, the same argument used in Remark \ref{R:l2estimate} shows that if $u\in\HH_{3,\theta}$
then $\chi(u(x,y,z))=|u(x,y,z)-\a_{+}|$ for $x\geq 0$ and $\chi(u(x,y,z))=|u(x,y,z)-\a_{-}|$ for $x< 0$. Hence,
 $|u(x,y,z)-z_{0}(x)|^{2}\leq 2\chi(u(x,y.z))^{2}+2\chi(z_{0}(x))^{2}$ on $\PP_{\theta}$. By (\ref{BGS}), there exists a constant $\tilde C>0$ such that if
$\|u\|_{L^{\infty}}\leq R$, then $|u(x,y)-z_{0}(x)|^{2}\leq \tilde CW(v(x,y))+2\chi(z_{0}(x))^{2}$.
 We derive that if $r>0$ and $u\in\HH_{3,\theta}$ is such that $\|u\|_{L^{\infty}}\leq R$ then
\begin{equation}\label{eq:stima3b}
\|u-z_{0}\|^{2}_{L^{2}(T_{r})}\leq 
\tilde C\varphi_{3}(u)+\tilde C (\tan(\theta) r^{2}m_{1}+rm_{2})+2\tan(\theta) r^{2}\int_{R}\chi(z_{0})^{2}\, dx.
\end{equation}
 }\end{Remark}

To obtain the existence of a minimum of $\F_{3}$ on $\HH_{3,\theta}$ it is now sufficient to apply the direct method of the Calculus of Variation.

\begin{Lemma}\label{L:existence} 
For every $\theta\in(0,\frac\pi4]$ there exists $u_{\theta}\in\MM_{3,\theta}$.
\end{Lemma}
\Proof Let $(u_{n})$ be a minimizing sequence for $\F_{3}$ in  $\HH_{3,\theta}$. Arguing as in Remark \ref{R:trapped}, we can assume that $\|u_{n}\|_{L^{\infty}(\PP_{\theta})}\leq R$
and so, by (\ref{eq:stima3a}) and (\ref{eq:stima3b}) we have that, for any $r>0$, $(u_{n}-z_{0})$ is a bounded sequence on $H^{1}(T_{r} )^{2}$. 
Then, by a diagonal argument, we obtain that there exist  $u\in z_{0}+\cap_{r>0}H^{1}(T_{r})$ and a subsequence of $(u_{n})$, still denoted $(u_{n})$, such that $u_{n}-u\to 0$ weakly in $H^{1}(T_{r})^{2}$ for any $r>0$ and for a.e. in $\PP_{\theta}$. 
 By the pointwise convergence $u$ enjoys the same symmetries of the functions $u_{n}$ and since $u\in z_{0}+\cap_{r>0}H^{1}(T_{r})^{2}$ we have indeed  that $u\in\HH_{3,\theta}$. The Lemma then follows by semicontinuity of $\F_{3}$. 
\QED

If $\psi\in C_{0}^{\infty}(\R^{3},\R^{2})$ verifies $\psi(-x,y,z)=\hat\psi(x,y,z)$ and $\psi(x,-y,z)=\tilde\psi(x,y,z)$ then $\F_{3}(u+\psi)\geq\F_{3}(u)$ for every $u\in\MM_{3,\theta}$. This property and the symmetry of the problem imply that every  $u\in\MM_{3,\theta}$ is a weak solution on $\PP_{\theta}$ of the system $-\Delta u+\nabla W(u)=0$ satisfying Neumann boundary condition on $\partial \PP_{\theta}$.
\begin{Lemma}\label{L:principiovariazionale} If $u$ is a minimum
of $\F_{3}$ on 
$\HH_{3,\theta}$ we have 
\[\int_{P_{\theta}}\nabla u\cdot\nabla\psi+\nabla W(u)\cdot\psi\, dx\, dy\, dz= 0\hbox{ for all }\psi\in C^{\infty}_{0}(\R^{3},\R^{2})\]
where we denote $\nabla u\cdot\nabla\psi=\sum_{i=1}^{2}\nabla u_{i}\cdot\nabla\psi_{i}$.
\end{Lemma}
\Proof
Given $h:\R^{3}\to\R$ we denote  
$\Lambda_{1}h(x,y,z)=\frac 1{2} (h(x,y,z)-h(-x,y,z))$ and $\Lambda_{2}h(x,y,z)=\frac 1{2} (h(x,y,z)-h(x,-y,z))$. We have that
$\Lambda_{i}h$ is odd with respect to the $i$-th variable while $(1-\Lambda_{i})h$ is even with respect to the $i$-th variable. 
Moreover,
$\Lambda_{i}\Lambda_{i}h=\Lambda_{i}h$, $\Lambda_{i}\Lambda_{j}h=\Lambda_{j}\Lambda_{i}h$ for $i\not=j$ and 
if $h$ is derivable then we have
$\partial_{i}\Lambda_{i}h=(1-\Lambda_{i})\partial_{i}h$ and  $\partial_{i}\Lambda_{j}h=\Lambda_{j}\partial_{i}h$, for $i\not=j$.\\
Finally, given $u:\R^{3}\to\R^{2}$, setting $u^{*}=((1-\Lambda_{2})\Lambda_{1}u_{1}, (1-\Lambda_{1})\Lambda_{2}u_{2})$,  note that $u\equiv u^{*}$ if and only if $u(-x,y,z)=\hat u(x,y,z)$ and $u(x,-y,z)=\tilde u(x,y,z)$. \smallskip

\noindent Given $u\in\MM_{3,\theta}$ and $\psi\in C^{\infty}_{0}(\R^{3},\R^{2})$,  we have $\psi^{*}\equiv(\psi^{*})^{*} $, so that,  $u+\psi^{*}\in\HH_{3,\theta}$. Then $\F_{3}(u+t\psi^{*})\geq\F_{3}(u)$
for any $t>0$
 and hence
\begin{equation}\label{eq:finalsyma}
\F_{3}( u+t\psi)-\F_{3}( u)\geq \F_{3}( u+t\psi)-\F_{3}( u+t\psi^{*}).
\end{equation}
Note that 
\begin{align*}
\int_{\PP_{\theta}}&|\nabla (u+t\psi)|^{2}-|\nabla (u+t\psi^{*})|^{2}\, dx dy dz\\
&=\int_{\PP_{\theta}}2t\, \nabla u\cdot\nabla(\psi-\psi^{*})+t^{2}\, (2\nabla\psi^{*}\cdot\nabla(\psi-\psi^{*})+|\nabla (\psi-\psi^{*})|^{2})\, dx dy dz.
\end{align*}
Since $u=((1-\Lambda_{2})\Lambda_{1}u_{1}, (1-\Lambda_{1})\Lambda_{2}u_{2})$, $\psi^{*}=((1-\Lambda_{2})\Lambda_{1}\psi^{*}_{1}, (1-\Lambda_{1})\Lambda_{2}\psi^{*}_{2})$ and $\psi-\psi^{*}=((1-\Lambda_{1})\psi_{1}+\Lambda_{2}\Lambda_{1}\psi_{1}, (1-\Lambda_{2})\psi_{2}+\Lambda_{1}\Lambda_{2}\psi_{2})$, using the above listed properties of the operators $\Lambda_{i}$, it can be shown that $\nabla u\cdot \nabla(\psi-\psi^{*})$ and $\nabla \psi^{*}\cdot \nabla(\psi-\psi^{*})$  are the sums of  functions which are odd with respect to one of the variables $x$ or $y$. By the symmetry of $\PP_{\theta}$ with respect to the plane $x=0$ and $y=0$, we obtain that $\int_{\PP_{\theta}}\nabla u\cdot \nabla(\psi-\psi^{*})dxdydz=\int_{\PP_{\theta}}\nabla\psi^{*}\cdot \nabla(\psi-\psi^{*})dxdydz=0$ from which 
\begin{equation*}\label{eq:finalsymb}
\int_{\PP_{\theta}}|\nabla (u+t\psi)|^{2}-|\nabla (u+t\psi^{*})|^{2}\, dx\, dy\, dz=\int_{\PP_{\theta}}t^{2}|\nabla (\psi-\psi^{*})|^{2}\, dx\, dy\, dz
\end{equation*}
and hence, by (\ref{eq:finalsyma}), that
\begin{equation}\label{eq:finalsyma*}
\F_{3}( u+t\psi)-\F_{3}( u)\geq \int_{\PP_{\theta}}W(u+t\psi)-W(u+t\psi^{*})\,dxdydz.
\end{equation}
We now note that, since $W(\xi_{1},\xi_{2})$ is even in both the variables,
the symmetries of $u$ imply that $(\partial_{1}W)(u(x,y,z))$ is odd in $x$ and even in $y$, while, $(\partial_{2}W)(u(x,y,z))$
is even in $x$ and odd in $y$.
In other worlds we have 
\[(\nabla W)(u(x,y,z))=((1-\Lambda_{2})\Lambda_{1}\partial_{1}W(u(x,u,z)), (1-\Lambda_{1})\Lambda_{2}\partial_{2}W(u(x,u,z)).\]
We then derive that
$\nabla W(u)\cdot(\psi-\psi^{*})=(1-\Lambda_{2})\Lambda_{1}\partial_{1}W(u)[(1-\Lambda_{1})\psi_{1}+\Lambda_{2}\Lambda_{1}\psi_{1}]+
(1-\Lambda_{1})\Lambda_{2}\partial_{2}W(u(x,u,z))[(1-\Lambda_{2})\psi_{2}+\Lambda_{1}\Lambda_{2}\psi_{2}]$, from which, again by symmetry, we deduce that
\begin{equation}\label{eq:finalsymc}
\int_{\PP_{\theta}}\nabla W( u)\cdot (\psi-\psi^{*})\, dx\, dy\, dz=0.
\end{equation}
By (\ref{eq:finalsyma*}) and (\ref{eq:finalsymc}) we conclude 
\begin{align*}
\int_{\PP_{\theta}}\nabla  u\cdot\nabla\psi+\nabla W( u)\cdot\psi\, &dx\, dy\, dz=\lim_{t\to 0^{+}}\tfrac{1}{t}(\F_{3}( u+t\psi)-\F_{3}( u))\\
&\geq
\lim_{t\to 0^{+}}\int_{\PP_{\theta}}\tfrac{W( u+t\psi)-W( u)}{t}+
\tfrac{W( u)-W( u+t\psi^{*})}{t}\, dx\, dy\, dz\\
&=\int_{\PP_{\theta}}\nabla W( u)\cdot (\psi-\psi^{*})\, dx\, dy\, dz=0.
\end{align*}
Since $\psi$ was arbitrary, the above argument shows that actually $\int_{\PP_{\theta}}\nabla  u\cdot\nabla\psi+\nabla W( u)\cdot \psi\, dx\, dy\, dz= 0$ for all $\psi\in C^{\infty}_{0}(\R^{3})$ and the Lemma follows.\QED

\noindent \begin{Remark}\label{R:3sol}{\rm By Lemma \ref{L:existence} and Lemma \ref{L:principiovariazionale}, using classical arguments, we obtain  that for every $\theta\in(0,\frac\pi4]$ there exists a a minimum $u_{\theta}$ of $\F_{3}$ on $\HH_{3,\theta}$ which is a classical solution  of the problem 
 \begin{equation*}
\begin{cases} -\Delta v(x,y,z)+\nabla W(v(x,y,z))=0,&
(x,y,z)\in \PP_{\theta},\cr 
v(x,y,z)=-v(-x,-y,z),& (x,y,z)\in \PP_{\theta}\cr
\partial_{\nu}v(x,y,z)=\bf 0 & (x,y,z)\in \partial \PP_{\theta}\cr
{\displaystyle\lim_{x\to\pm\infty}|v(x,y,z)-{\bf a}_{\pm}|=0}, 
	&\hbox{for }(y,z)\,/\, z>0\hbox{ and }-\bar z\leq y\leq \bar z
\end{cases}\eqno(P_{3,\theta}).
 \end{equation*} 
  Since $u_{\theta}\in\MM_{3,\theta}$ we have $\|u_{\theta}\|_{L^{\infty}}\leq R$  and by Schauder estimates there exists of a constant $C>0$ ( $C$ is independent of $\theta$ since, by recursively reflecting the solution $u_{\theta}$ with respect to the faces of $\PP_{\theta}$, see below, we can always think at $u_{\theta}$ as the restriction to $\PP_{\theta}$ of a solution $u$ of $-\Delta v(x,y,z)+\nabla W(v(x,y,z))=0$ on $z>0$ with $\|u\|_{L^{\infty}}\leq R$) such that
\begin{equation}\label{eq:c2stime}
\|u_{\theta}\|_{C^{2}(\PP_{\theta})}\leq C\quad \hbox{for any }\theta\in (0,\pi/4].
\end{equation}
The estimate (\ref{eq:c2stime}) and the fact that $\F_{3}(u_{\theta})<+\infty$ allow to show, via an argument analogous to the one used in proving Remark \ref{R:unif}, that $\lim_{x\to\pm\infty}|v(x,y,z)-{\bf a}_{\pm}|=0$, uniformly on $\{(y,z)\,|\, z>0\hbox{ and }-\bar z\leq y\leq \bar z\}$.}
\end{Remark} 

We now study the asymptotic behaviour of $u_{\theta}\in\MM_{3,\theta}$ as $z\to +\infty$.

\begin{Remark}\label{R:estimate3} {\rm  
The functionals $\F_{3,I}$ are iteratively defined starting from the lower dimensional functionals $\F_{1}$ and $\F_{2}$ and some of the estimates we gave in previous sections have the analogous in the 3 dimensional case. 
Given $u\in \HH_{3,\theta}$ we have
\begin{itemize}
\item[$(i)$] for all $(z_{1},z_{2})\subset\R^{+}$ there results
$$
\|u(\cdot,\cdot,z_{2})-u(\cdot,\cdot, z_{1})\|^{2}_{L^{2}(S_{\bar z_{1}})^{2}}
	\leq 2\F_{3,(z_{1},z_{2})}(u) |z_{2}-z_{1}|.
$$

	\item[$(ii)$] If $\F_{2,\bar z}(u(\cdot,\cdot, z))- m_{2,\bar z}\geq \nu>0$ for almost every $z\in(\sigma,\tau)\subset\R^{+}$, then
\begin{align*}
\varphi_{3,(\sigma,\tau)}(u)&
\geq \frac{1}{2(\tau-\sigma)}\|u(\cdot,\cdot, \tau)-u(\cdot,\cdot, \sigma)\|^{2}_{L^{2}(S_{\bar \sigma})^{2}}+\nu(\tau-\sigma)\\
&\geq
\sqrt{2\nu}\|u(\cdot,\cdot,\tau)-u(\cdot,\cdot,\sigma)\|_{L^{2}(S_{\bar \sigma})^{2}}.
\end{align*}
\end{itemize}}
\end{Remark}

\noindent The estimates given in Remark \ref{R:estimate3}  allow us to obtain  the following result

\begin{Lemma}\label{L:asymptot3}
If $u\in\HH_{3,\theta}$ and $\F_{3}(u)<+\infty$, then, there exists $q\in\MM_1^{*}$ such that $\dist_{L^{2}(S_{\bar z})^{2}}(u(\cdot,\cdot,z),\MM_{2,q})\to 0$ as $z\to +\infty$.
\end{Lemma}
\Proof 
Since $\F_{3}(u)<+\infty$  
there exists $\bar\ell+1<z_{n}\to+\infty$ such that
$\F_{2,\bar z_{n}}(u(\cdot,\cdot,z_{n}))-m_{2,\bar z_{n}}\to 0$. Then, 
by Lemma \ref{L:concentrazione2dL}  there exists $q\in\MM_1^{*}$ such that, up to a subsequence, $\|u(\cdot,y,z_{n})-q\|_{L^{2}}\le d_{0}$ for all $y\in(\bar\ell,\bar z_{n})$ and $n\in\N$, and then, by Lemma \ref{L:compattezzaperiodiche}, we obtain that 
$\dist_{L^{2}(S_{\bar z_{n}})^{2}}(u(\cdot,\cdot,z_{n}),\MM_{2,q})\to 0$ as $n\to +\infty$. \\
Assume now by contradiction that
$\dist_{L^{2}(S_{\bar z})^{2}}(u(\cdot,\cdot,z),\MM_{2,q})\not\to 0$ as $z\to+\infty$. 
Then  there exist $\delta>0$ and 
$(\zeta_{n})\subset (0,+\infty)$ such that $\dist_{L^{2}(S_{\bar \zeta_{n}})^{2}}(u(\cdot,\cdot,\zeta_{n}),\MM_{2,q})\geq \delta$ for any $n\in\N$.
We can assume, extracting subsequences if needed, that $z_{n}< \zeta_{n}<z_{n+1}$ for any $n\in\N$.
Since $\dist_{L^{2}(S_{\bar \zeta_{n}})^{2}}(u(\cdot,\cdot,z_{n+1}),\MM_{2,q})\to 0$ as $n\to +\infty$, by Remark \ref{R:estimate3}-$(i)$
there exist $\bar n\in\N$, a sequence of intervals $(\s_{n},\tau_{n})\subset (\zeta_{n},z_{n+1})$ and   $\rho\in (0,d_{0})$ for which
\begin{itemize}
\item[(i)] $\|u(\cdot,\cdot,\tau_{n})-u(\cdot,\cdot,\sigma_{n})\|_{L^{2}(S_{\zeta_{n}})^{2}}=\rho$, for all $n\geq\bar n$
\item[(ii)] $2\rho\geq\dist_{L^{2}(S_{\zeta_{n}})^{2}}(u(\cdot,\cdot,z),\MM_{2,q})\geq \rho$, for every $z\in (\s_{n},\tau_{n})$, $n\geq\bar n$.
\end{itemize}
Since, by Remark \ref{R:distM}, $\dist_{L^{2}(S_{\zeta_{n}})^{2}}(\MM_{2,q},\MM_{2,p})\geq 3d_{0}$ for all $p\not=q\in\MM_1^{*}$, by (ii) we get $\min_{p\in\MM_1^{*}}\dist_{L^{2}(S_{\zeta_{n}})^{2}}(u(\cdot,\cdot,z),\MM_{2,p})\geq \rho$ for every $z\in (\s_{n},\tau_{n})$ and hence, by Corollary \ref{C:mu0L2}, we recover that there exists $\mu_{\rho}>0$  such that, taking $\bar n$ larger if necessary, $\F_{2,\bar z}(u(\cdot,\cdot,z))-m_{2,\bar z}\geq \mu_{\rho}$ for all $z\in (\s_{n},\tau_{n})$ and $n\geq\bar n$.
Using now Remark \ref{R:estimate3}-(ii) we conclude 
$\varphi_{3,(\s_{n},\tau_{n})}(u)\geq\sqrt{2\mu}\tilde\delta>0$ for all $n\geq \bar n$ and so $\varphi_{3}(u)=+\infty$, a contradiction.\QED

The Schauder estimates (\ref{eq:c2stime}) allow to improve the result in Lemma \ref{L:asymptot3} for the functions $u\in\MM_{3,\theta}$.
\begin{Lemma}\label{L:asymptot4}
If $u\in\MM_{3,\theta}$ then  there exists $q\in\MM_1^{*}$ such that 
$$
\dist_{L^{\infty}(S_{\bar z})^{2}}(u_{\theta}(\cdot,\cdot,z),\MM_{2,q})\to0\,\hbox{ 
as $z\to+\infty$}.
$$
\end{Lemma}
\Proof By Lemma \ref{L:asymptot3} there exists $q\in\MM_1^{*}$ such that $\dist_{L^{2}(S_{\bar z})^{2}}(u(\cdot,\cdot,z),\MM_{2,q})\to 0$ as $z\to +\infty$. Assume by contradiction that there exists $\rho_{0}>0$ and  $(x_{n},y_{n},z_{n})\in\PP_{\theta}$ with $z_{n}\to +\infty$ such that fixed any $v\in\MM_{2,q}$ there results $
|u(x_{n},y_{n},z_{n})-v(x_{n},y_{n})|\geq 4\rho_{0}$.
 By (\ref{eq:c2stime}) we obtain that there exists $r_{0}>0$ such that if $n\in\N$ and $(x,y,z)\in\PP_{\theta}$ is such that $|x-x_{n}|, |y-y_{n}|, |z-z_{n}|\leq 2r_{0}$ then
$|u(x,y,z)-v(x,y)|\geq 2\rho_{0}$. Then $\|u(\cdot,\cdot,z_{n})-v(\cdot,\cdot)\|_{L^{2}(S_{\bar z_{n}})^{2}}\geq 8\rho_{0}r_{0}$ for every $n\in\N$ and
$v\in\MM_{2,q}$, a contradiction which proves the Lemma.\QED

To prove Theorem \ref{T:main} we now fix $j\in\N$, $j\geq 2$ and $\theta_{j}={\frac\pi{2j}}$. Letting  $u_{j}\equiv u_{\theta_{j}}$ be given by Lemma \ref{L:existence}, we construct  an entire solution $v_{{j}}:\R^{3}\to \R^{2}$ to (\ref{eq:P3}) satisfying the properties stated in the Theorem \ref{T:main} recursively reflecting $u_{j}$ with respect to the faces of $\PP_{\theta_{j}}$.  
More precisely, for $j\in\N$, $j\geq 2$, we consider  the rotation matrix
$$
A_{j}=\left\|\begin{array}{ccc} 1&0&0\\ 0 & \cos\frac{\pi}{j} & \sin\frac{\pi}{j} \\ 0&-\sin\frac{\pi}{j} & \cos\frac{\pi}{j} \\\end{array}\right\|
$$
Setting $\PP_{k,j}=A_{j}^{k}\PP_{\frac\pi{2j}}$, for every $k=0,...,2j-1$, we have
$\R^{3}=\cup_{k=0}^{2j-1}\PP_{k,j}$
and that if $k_{1}\not=k_{2}$ then $\hbox{int}(\PP_{k_{1},j})\cap \hbox{int}(\PP_{k_{2},j})=\emptyset$.	
	
\noindent We consider the function $\bar u_{j}(x,y,z)=u_{j}(x,-y,z)$, the reflected of $u_{j}$ with respect to $y=0$. 
Noting that $A^{-k}_{j}\PP_{k,j}= \PP_{\frac{\pi}{2j}}$, for $(x,y,z)\in P_{k,j}$, $k=0,...,2j-1$, we define
$$
v_{j}(x,y,z)=u_{j}(A^{-k}_{j}(x,(-1)^{k}y,z))
=\begin{cases}u_{j}(A^{-k}_{j}(x,y,z)),&\hbox{if $k$ is even,}\cr \bar u_{j}(A^{-k}_{j}(x,y,z)),&\hbox{if $k$ is odd.}\end{cases}
$$
We recognize that $v_{j}|_{\PP_{1,j}}$ is the reflection of $v_{j}|_{\PP_{0,j}}=u_{j,q}$ w.r.t. the plane which separates $\PP_{0,j}=\PP_{\frac\pi{2j}}$ from $\PP_{1,j}$. Recursively, for every $k\in\{ 1,\ldots, 2j-1\}$, we have that $v_{j}|_{\PP_{k,j}}$ is the reflection of $v_{j}|_{\PP_{k-1,j}}$ w.r.t. the plane separating $\PP_{k-1,j}$ from $\PP_{k,j}$. Note that $v_{j}\in H^{1}_{loc}(\R^{3})$  (see \cite{[Brezis]} Lemma IX.2.).
Moreover, note that  if $\psi\in C^{\infty}_{0}(\R^{3},\R^{2})$ and $k\in\{1,\ldots,2j-1\}$ then, trivially,
$\psi\circ A_{j}^{k}\in C^{\infty}_{0}(\R^{3},\R^{2})$ and so by Lemma \ref{L:principiovariazionale},  
we obtain
\begin{align*}
\int_{\PP_{k,j}}\nabla v_{j}(x,y,z)&\cdot\nabla\psi(x,y,z)+\nabla W(v_{j}(x,y,z))\cdot \psi(x,y,z)\, dx\, dy\, dz\\&=
\int_{\PP_{\frac\pi{2j}}}\nabla u_{j}(x,(-1)^{k}y,z)\cdot\nabla\psi\circ A_{j}^{k}(x,(-1)^{k}y,z)\\&\quad
+\nabla W(u_{j}(x,(-1)^{k}y,z))\cdot \psi\circ A_{j}^{k}(x,(-1)^{k}y,z)\, dx\, dy\, dz=0.
\end{align*}
We conclude that for any $\psi\in C^{\infty}_{0}(\R^{3},\R^{2})$, we have
$$
\int_{\R^{3}}\nabla v_{j}\cdot\nabla\psi+\nabla W(v_{j})\psi\, dx\, dy=
\sum_{k=0}^{2j-1}\int_{\PP_{k,j}}\nabla v_{j}\cdot\nabla\psi+\nabla W(v_{j})\psi\, dx\, dy=0,
$$
i.e., $v_{j}$ is a weak and so, by standard bootstrap arguments, a classical solution of $-\Delta v+\nabla W(v)=0$ on $\R^{3}$. We finally recall that, by Lemma \ref{L:asymptot3}, there exists $ q_{j}\in\MM_1^{*}$ 
such that
$\dist_{L^{\infty}(S_{\tan(\frac{\pi}{j})z})^{2}}(u_{j}(\cdot,\cdot,z),\MM_{2,q_{j}})\to0$ 
as $z\to+\infty$ and gathering all the above results, denoting 
$\breve{v}_{j}(x,\rho,\theta)=v_{j}(x,\rho\cos\theta,\rho\sin\theta)$, we conclude that
$v_{j}$ satisfies the properties:
\begin{itemize}
\item[(i)]  $\dist_{L^{\infty}(S_{\tan(\frac{\pi}{j})z})^{2}}(v_{j}(\cdot,\cdot,z),\MM_{2,q_{j}})\to 0$ as $z\to +\infty$.
\item[(ii)] $\breve{v}_{j}$ is periodic in $\theta$ with period $\frac{2\pi}{j}$
\item[(iii)]  ${\displaystyle\lim_{\rho\to+\infty}}\breve{v}_{j}(x,\rho,\tfrac\pi 2+\tfrac\pi{j}(\tfrac 12+k))=\begin{cases}q_{j}(x)&\hbox{if }k\hbox{ is odd}\\
\bar q_{j}(x)&\hbox{if }k\hbox{ is even}\end{cases}$ uniformly
for $x\in\R$.
\end{itemize}

\end{document}